\documentclass[]{aiaa-tc}

 \title{High Order Flux Reconstruction on Stretched and Warped Meshes}

 \author{
  Will Trojak%
    \thanks{PhD Candidate, Department of Engineering, University of Cambridge, AIAA Student Member.}
   , Rob Watson%
   	\thanks{Research Fellow, Department of Engineering, University of Cambridge, AIAA Member}
  \ and Paul G. Tucker%
  	\thanks{Professor, Department of Engineering, University of Cambridge, AIAA Assoc. Fellow}  
  \thanksibid{0}\\
  {\normalsize\itshape
   Department of Engineering, University of Cambridge, Cambridge, UK, CB2 1PZ}\\
 }

 \AIAApapernumber{YEAR-NUMBER}
 \AIAAconference{Conference Name, Date, and Location}
 \AIAAcopyright{\AIAAcopyrightD{YEAR}}

 \newcommand{\etal}{\emph{et~al.}}

\graphicspath{ {Figs/} }

\usepackage{graphicx}
\usepackage{graphics}
\usepackage{upgreek}
\usepackage{float}
\usepackage{epstopdf}
\usepackage{amsmath}
\usepackage{wrapfig}
\usepackage{amssymb}
\usepackage{pdflscape}
\usepackage{listings}
\usepackage{cite}
\usepackage{subcaption}
\usepackage{multirow}
\usepackage[percent]{overpic}
\usepackage{multicol}
\usepackage{soul}
\usepackage{xcolor}

\definecolor{c8_1}{rgb}{0.698039215686275,0.0941176470588235,0.168627450980392}
\definecolor{c8_2}{rgb}{0.839215686274510,0.376470588235294,0.301960784313725}
\definecolor{c8_3}{rgb}{0.956862745098039,0.647058823529412,0.509803921568627}
\definecolor{c8_4}{rgb}{0.992156862745098,0.858823529411765,0.780392156862745}
\definecolor{c8_5}{rgb}{0.819607843137255,0.898039215686275,0.941176470588235}
\definecolor{c8_6}{rgb}{0.572549019607843,0.772549019607843,0.870588235294118}
\definecolor{c8_7}{rgb}{0.262745098039216,0.576470588235294,0.764705882352941}
\definecolor{c8_8}{rgb}{0.129411764705882,0.400000000000000,0.674509803921569}

\begin{document}

\maketitle

\begin{abstract}
	High-order CFD is gathering a broadening interest as a future industrial tool, with one such approach being Flux Reconstruction (FR). However, due to the need to mesh complex geometries if FR is to displace current, lower order methods, FR will likely have to be applied to stretched and warped meshes. Therefore, it is proposed that the analytical and numerical behaviour of FR on deformed meshes for both the 1D linear advection and the 2D Euler equations is investigated. The analytical foundation of this work is based on a modified von Neumann analysis for linearly deformed grids that is presented. The temporal stability limits for linear advection on such grids are also explored analytically and numerically, with CFL limits set out for several Runge-Kutta schemes, with the primary trend being that contracting mesh regions give rise to higher CFL limits whereas expansion leads to lower CFL limits. Lastly, the benchmarks of FR are compared to finite difference and finite volumes schemes, as are common in industry, with the comparison showing the increased wave propagating ability on warped and stretched meshes, and hence, FR;s increased resilience to mesh deformation. 
\end{abstract}

\section*{Nomenclature}
\begin{multicols}{2}
\begin{tabbing}
  XXXX \= \kill
  	\textit{Roman}\\
	$a$ \> convective velocity \\	  
	$c(k)$ \> phase velocity at wavenumber $k$ \\
	$\mathbf{C}_0$ \> centre cell FR matrix \\
	$\mathbf{C}_{-1}$ \> upwind cell FR matrix \\
	$\mathbf{D}$ \> first derivative matrix \\
	$f$ \> flux variable in physical domain \\
	$\hat{G}(\hat{k})$ \> computational filter kernel \\	
	$h_l \:\mathrm{\&}\: h_r$ \> left and right correction functions\\
	$J_n$ \> $n^{\mathrm{th}}$ cell Jacobian\\
	$k$ \> wavenumber \\
	$k_{nq}$ \> solution point Nyquist wavenumber, $(p+1)/\delta_j$\\
	$\hat{k}$ \> $k_{nq}$ normalised wavenumber, $[0,\pi]$ \\
	$l_n$ \> $n^{\mathrm{th}}$ Lagrange basis function \\
	$p$ \> solution polynomial order \\
	PPW \> points per wavenumber at given error level \\
	$\mathbf{Q}$ \> spatial scheme matrix \\
	$\mathbf{R}$ \> update matrix \\
	$u$ \> primitive in physical domain \\
	$\tilde{u}$ \> Fourier reconstructed field \\		
	
	\textit{Greek}\\
	$\gamma$ \> grid geometric expansion factor \\
	$\delta_j$ \> mesh spacing, $x_j-x_{j-1}$\\
	$\epsilon$ \> PPW error level\\	
	$\xi$ \> transformed spatial variable \\
	$\rho(\mathbf{A})$ \> spectral radius of $\mathbf{A}$ \\ 
	$\tau$ \> time step \\
	$\pmb{\Omega}$ \> solution domain \\
	$\pmb{\Omega}_n$ \> $n^{\mathrm{th}}$ solution sub-domain \\
	$\pmb{\Omega}_s$ \> standardised sub-domain\\
	
	\textit{Subscript}\\
	$\mathrm{\bullet}_l$ \> variable at left of cell\\
	$\mathrm{\bullet}_r$ \> variable at right of cell\\
	
	\textit{Superscript}\\
	$\mathrm{\bullet}^c$ \> common value at interface\\
	$\mathrm{\bullet}^T$ \> vector or matrix transpose \\
	$\mathrm{\bullet}^{\delta}$ \> discontinuous value\\
	$\hat{\mathrm{\bullet}}$ \> transformed variable \\
	$\overline{\mathrm{\bullet}}$ \> locally fitted polynomial of variable
\end{tabbing}
\end{multicols}

\section{Introduction}
\begin{figure}
		\centering
			\includegraphics[width=0.43\linewidth,trim= 31mm 82mm 31mm 73mm,clip=true]{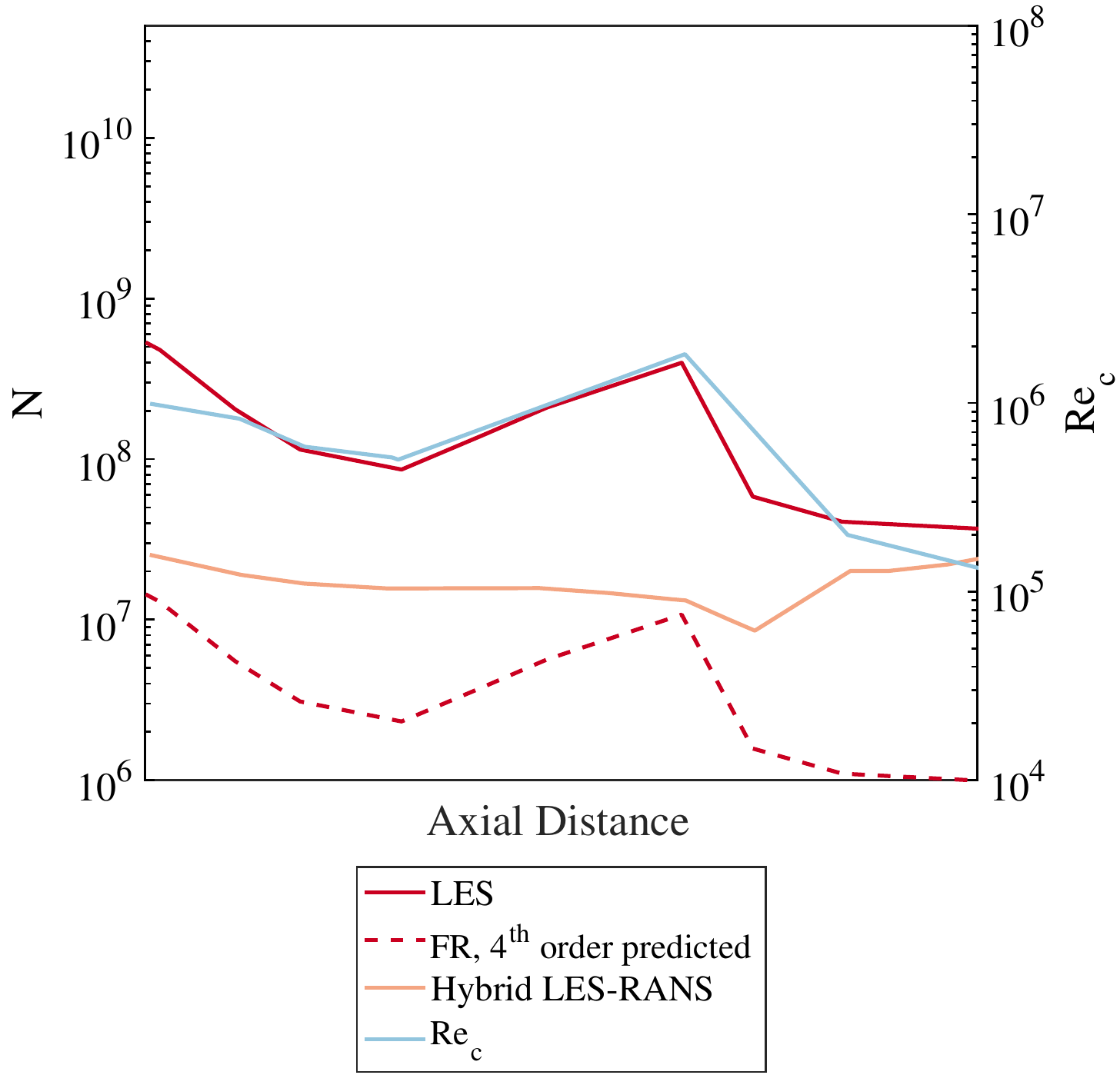}
			\put(-173.5,135){\includegraphics[width=0.310\linewidth, trim= 67mm 107mm 75mm 65mm,clip=true]{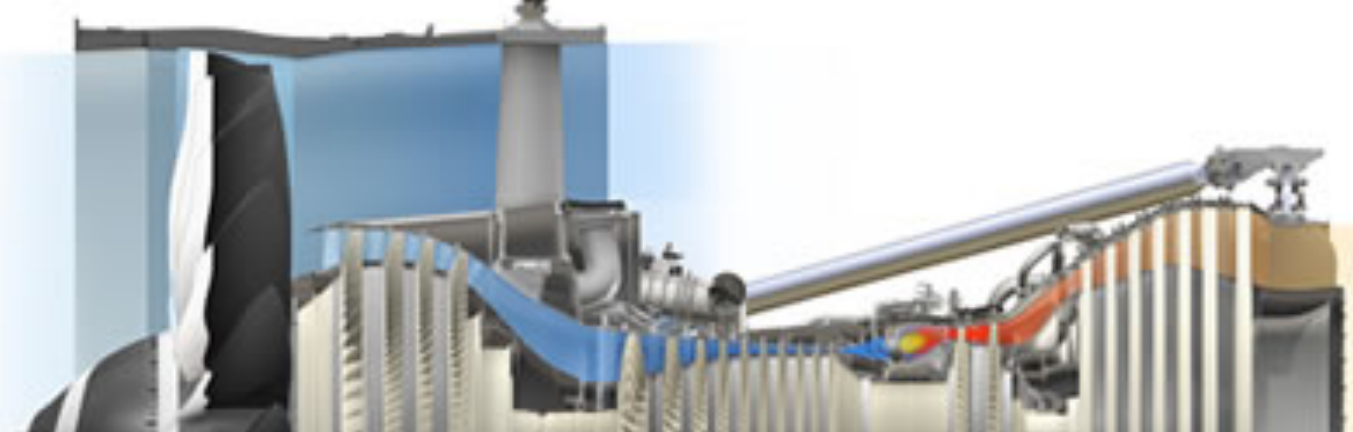}}
		\captionof{figure}{LES and hybrid RANS-LES mesh degrees of freedom requirement against axial distance through a gas turbine engine. $Re$  is also plotted, showing large variations in flow regimes occur. (Tucker et al.~\cite{Tucker2010}) This shows that with current methods LES is only routinely used for the low pressure turbines. Data presented later can give an estimate of the 3D mesh requirement for $4^{th}$ order Flux Reconstruction (FR).}
		\label{fig:LES-RANS}
	\end{figure}

The potential of Large Eddy Simulation (LES) has been understood for some time, however most current industrial LES and Reynolds-Averaged Navier Stokes (RANS) implementations make use of second order spatial schemes. This kind of lower spatial order for LES can be prohibitively expensive in some flow regimes, with cost scaling, from Piomelli~\cite{Piomelli2008a}, with $\sim Re^{2.4}$ for the innermost section of the boundary layer. Figure~\ref{fig:LES-RANS} aims to show that for a gas turbine engine conventional LES has a much larger overhead than hybrid RANS-LES, with full LES only really being currently feasible in the lower pressure turbine stages. From metrics shown later, using FR could make conventional LES far cheaper. The high cost of wall-resolved LES originates from all but 10-20\% of the vortical motions being directly simulated in both space and time, Tucker~\cite{Tucker2014}. The unresolved fraction of the vortical motion must, therefore, be accounted for using a sub-grid scale model. As LES will make use of finite discretisations on the underlying equation (be that through a Finite Volume (FV), Discontinuous Galerkin (DG) or another method) the solution  will incur truncation and aliasing error. Therefore, for the sub-grid scale model to correctly influence the flow, the sub-grid-scale error must be minimised. This area was explored by Chow \& Moin~\cite{Chow2003} and Ghosal~\cite{Ghosal1996}, who found that for a second order central FD LES the filter width was needed to be four times the grid spacing. They also found that for an eighth order central FD LES a filter width of twice the grid spacing was required. This effect occurs because for second order schemes the numerical error and the sub-grid-scale forcing will have the same order scaling with respect to the grid spacing, whereas a move to higher order means the differencing error order increases and so the filter-scale can decrease. Regardless, higher orders allow for coarser meshes to capture the same flow regime, and can reduce the high computational barrier to the use of LES. Hence high order methods have become of increasing interest for industry. For example, studying the abstracts of AIAA SciTech Conferences on Aerospace Sciences shows  an increase from 62 to 99 papers that used or studied high order methods between 2014 and 2015. 

Huynh \cite{Huynh2007} introduced a high order scheme which he termed Flux Reconstruction (FR), and this has since evolved into a wider family of schemes explored by Vincent~\etal~\cite{Vincent2011} using energy stability arguments to define a wide family of correction functions. Extensions have been made to handle advection-diffusion~\cite{Huynh2009,Williams2014a}, simplex and hypercube elements\cite{Huynh2011,Williams2014,Sheshadri2016}, and various optimisations of the collocation points and correction functions~\cite{Vincent2015,Asthana2014,Asthana2015a}. The quantitative evaluation performed analytically and numerically has largely been confined to canonical meshes and the scenarios explored, although of great importance, are often of less direct relevance to CFD practitioners. This is due to their contrived nature leading to well presented problems, that are useful as benchmarks, but do not always show the full picture.  In essence, it is common for the actual problems encountered by engineers to be highly complex, making meshes with low levels of skew or low level inter-element expansion very difficult to produce. Furthermore, FR is rooted in finite element methods and consequently is an unstructured approach that fundamentally enables far more irregularities within the mesh compared to structured approaches. This highlights the importance of accurate characterisation of the scheme on irregular meshes.

Previous investigations into the analytical behaviour of numerical methods on irregular grids has been confined to schemes such as finite difference(FD), where the mesh transformation can be more easily applied. For example Chung~\&~Tucker~\cite{Chung2003} in which the effect of hyperbolically transformed grids was investigated for FD and compact FD schemes, clearly demonstrated the added dissipation that transformation can cause. More recently You~\etal~\cite{You2006} demonstrated a more thorough approach, detailing the exact higher order terms which lead to inaccuracy on hyperbolically skewed grids. Within the field of FR, some numerical results have been presented pertaining to the accuracy of FR applied to curvelinear grids, Mengaldo~\etal~\cite{Mengaldo2016}, also showing some interesting parallels between FR and Discontinuous Galerkin methods. This did, however, limit the quantification to numerical experiments, which, although very powerful, do not allow the bedrock of the scheme to be exposed. Further investigation into the effect of curvilinear grids on the Jacobian was presented by Kopriva~\cite{Kopriva2006}. Lastly, FR coupled to \emph{r}-type mesh adpatation was performed by Sheshadri~\etal~\cite{Sheshadri2015} for the purposes of shock capture using a 'divide or conquer' method that showed reasonable performance. These works have laid a fine basis for the development of study into the effect of mesh transformation on finite element schemes and in this case for FR. 

The aims of this paper are broken down into several sections. An analytical framework is constructed using the established von Neumann analysis, however generality of the grid will be maintained such that non-uniform grids may be investigated. The subsequent analytical investigations are concerned with the dispersion and dissipation of various orders of FR for non-uniform 1D grids, and with moving on to couple the spatial scheme to a temporal scheme to establish theoretical CFL limits for various grids. A numerical methodology is proposed such that the analytical results may validated, also enabling evaluation of the fully discretised spatial-temporal scheme. Lastly the numerical investigation is taken further, into 2D, using the Euler equations to understand the behaviour of more relevant flows on meshes that have been artificially degraded. This test case can also enable comparison to be made to a more prevalent industrial FV scheme.

\section{Flux Reconstruction}
	Flux Reconstruction ~\cite{Castonguay2012,Huynh2007} (FR) applied to the linear advection equation will form the basis of the initial investigation to be carried out, and, for the readers' convenience, an overview of the scheme is presented here. However, for a more detailed understanding the reader should consult Castonguay \cite{Castonguay2012} or Huynh \cite{Huynh2007}. This 1D scheme can be readily converted to two dimensions (carried out later) and three dimensions for quadrilaterals and hexahedrals respectively. First, let us consider the one dimensional advection equation:
	\begin{equation}
		\frac{\partial u}{\partial t} + \frac{\partial f}{\partial x} = 0
	\end{equation}
	
The FR method is related to the Discontinuous Galerkin (DG) method~\cite{Reed1973} and utilises the same subdivision of the domain into discontinuous sub-domains:
	\begin{equation}\label{eq:Domain}
		\pmb{\Omega}  = \bigcup_{n=1}^{N}{\pmb{\Omega}_n}
	\end{equation}

\begin{figure}
		\centering
		\begin{subfigure}[b]{0.35\linewidth}
			\centering
			\includegraphics[width=\linewidth,trim= 0mm 0mm 0mm 0mm,clip=true]{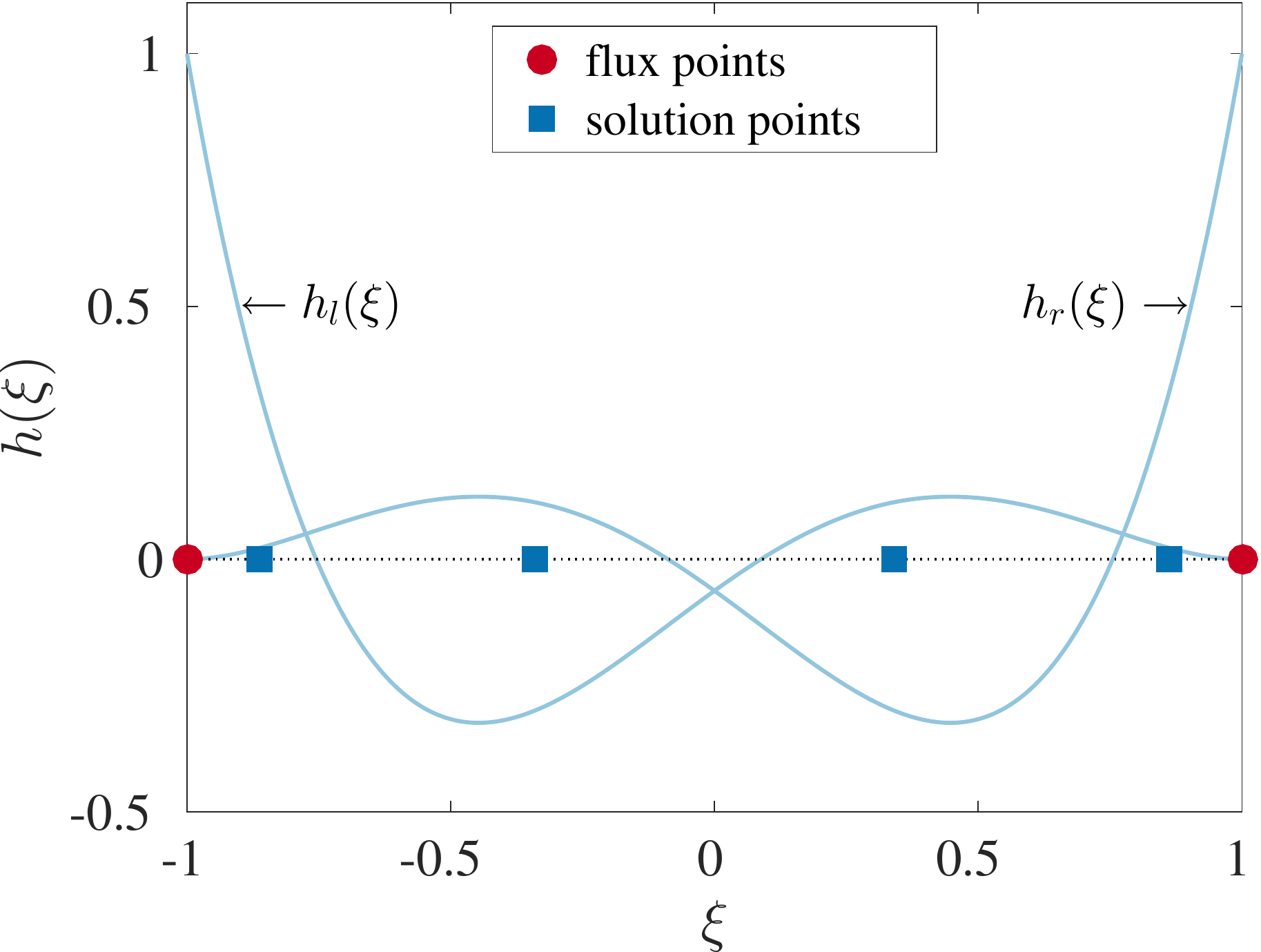}
			\caption{Flux and solution point layout for $p = 3$ in $\pmb{\Omega}_s$, with corresponding left and right Huynh, $g_2$, correction functions~\cite{Huynh2007}.}
			\label{fig:1D_layout}
		\end{subfigure}
		~
		\begin{subfigure}[b]{0.5\linewidth}
			\centering
			\includegraphics[width=\linewidth]{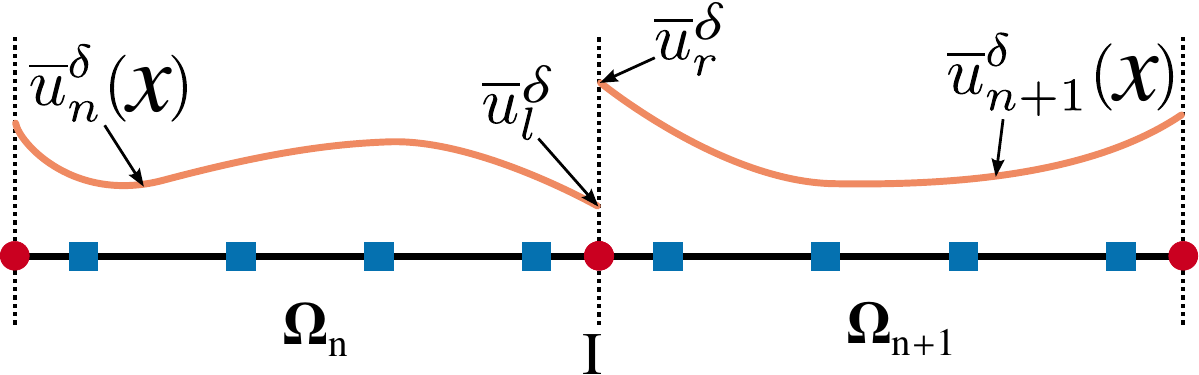}
			\caption{Discontinuous primitive polynomials and the interpolated values at the flux points for adjacent sub-domains.\\ \:\:\:}
			\label{fig:1D_interface_u}
		\end{subfigure}
		
		\caption{Point layout in $\pmb{\Omega}_s$ for $p=3$ and cell interface topology.}
		\label{fig:fr_interface}
	\end{figure}

	In the standardised sub-domain, $\pmb{\Omega}_s \in \mathbb{R}^d$, computational spatial variables are defined. When $d=1$,  $\pmb{\Omega}_s = [-1,1]$, using $\xi$ to denote the value taken. This computational space is discretised with $(p+1)^d$ solution points, and $2d(p+1)^{d-1}$ flux points are placed at the edges of the sub-domain. (Figure~\ref{fig:1D_layout} shows a 1D example). To transform from $\pmb{\Omega}_n \rightarrow \pmb{\Omega}_s$ the Jacobian $J_n$ is defined such that
	\begin{equation}
		\hat{u}^{\delta} = \hat{u}^{\delta}(\xi;t) = J_n u^{\delta}(x;t) 
	\end{equation}
	and the solution point mapping in the physical domain used throughout this investigation results in the Jacobian representing a linear mapping. With the domain set up, we now proceed with defining the steps to construct a continuous solution from discontinuous segments. The first stage is to define a local solution polynomial in $\pmb{\Omega}_s$ using Lagrange interpolation.
	\begin{align}
		l_k(\xi) &= \prod_{i=1,i\ne k}^{p+1}{\frac{\xi - \xi_i}{\xi_k - \xi_i}}\\
		\hat{u}^{\delta}(\xi;t) &= \sum_{i=1}^{p+1}{\hat{u}^{\delta}_il_i(\xi)}
	\end{align}	
	Repeating the interpolation for the discontinuous flux in $\pmb{\Omega}_s$:
	\begin{equation}
		\hat{f}^{\delta} = \hat{f}^{\delta}(\xi;t) = \sum_{i=1}^{p+1}{\hat{f}^{\delta}_il_i(\xi)}
	\end{equation}	
	Now using the Jacobian and the solution polynomials, the primitive and flux values can be approximated in the physical domain $\pmb{\Omega}_n$
	\begin{equation}
		\overline{u}^{\delta}(x;t) = \frac{\hat{u}^{\delta}(\xi;t)}{J_n}
	\end{equation}
	This distinction is made due to the potentially approximate nature of $J_n$. However for the case of linear transformations, as will be considered here, $J_n$ fully captures the spatial transformation and thus $\overline{u}^{\delta}$ will be a polynomial of order $p$. With a polynomial formed in the reference domain the value of the primitive at the flux points, see Fig.~\ref{fig:1D_layout}, can be defined as $\hat{u}^{\delta}_l = \hat{u}^{\delta}(-1)$ and $\hat{u}^{\delta}_r = \hat{u}^{\delta}(1)$, this can be repeated for the flux values.		
	Once the primitive values at the interface, $I$, have been interpolated, a common interface flux can be calculated, $f^{c}_I$, in the physical domain. For a general case this is done using a Riemann solver on the primitives at the interface, such as: Roe~\cite{Roe1981}; flux-vector splitting~\cite{vanLeer1982}; or HLL~\cite{Harten1983}. In order to get a spatially continuous solution over $\pmb{\Omega}$ the common interface flux must be incorporated into the solution. For FR this is done by using a correction function to propagate the corrected flux gradient into the $\pmb{\Omega}_n$. Figure~\ref{fig:1D_layout} shows the correction function  proposed by Huynh~\cite{Huynh2007} and, in general, the left and right correction functions are defined as $h_l(\xi)$ and $h_r(\xi)$. The procedure to apply the correction to transformed discontinuous flux in $\pmb{\Omega}_s$ is:
	\begin{equation}\label{eq:discrete_1D_full}
		\frac{\partial \hat{f}(\xi_i)}{\partial \xi} = \sum_{j=1}^{p+1}{\hat{f}_j^{\delta}\frac{dl_j(\xi_i)}{d\xi}} + (\hat{f}^{c}_l - \hat{f}^{\delta}_l)\frac{dh_l(\xi_i)}{d\xi} + (\hat{f}^{c}_r - \hat{f}^{\delta}_r)\frac{dh_r(\xi_i)}{d\xi}
	\end{equation}
	This then allows then transformed continuous equation can be written as:
	\begin{equation} \label{eq:discrete_1D_ad}
		\frac{\partial \hat{u}(\xi_i)}{\partial t} = - \frac{\partial \hat{f}(\xi_i)}{\partial \xi}
	\end{equation}
	and hence the solution can be advanced in time via some method of temporal integration.
	
\section{von Neumann Analysis}\label{sec:vn}
	The main analysis to be carried out is a von Neumann analysis that follows the work of Lele~\cite{Lele1992}, Hesthaven~\etal~\cite{Hesthaven2007}, Huynh~\cite{Huynh2007}, Vincent~\etal~\cite{Vincent2011}, and Asthana~\etal~\cite{Asthana2014}. The analysis shown here differs from that of previous work in that it does not assume that the grid is uniform and instead retains generality, the advantage being that the effect of grid stretching can be investigated. For the  analysis, consider the computational domain as before, $\pmb{\Omega} \in \mathbb{R}^1$ with $N-1$ sub-domains $\pmb{\Omega}_n$, with flux points of $\pmb{\Omega}_n$ located at $\mathbf{x}_j \: \forall \:\{j \in \mathbb{N},\: j \leqslant N\}$. The 1D linear advection  equation can be written as Eq.~(\ref{eq:1dadvection}).
	
	\begin{equation}\label{eq:1dadvection}
		\frac{\partial u}{\partial t} + a\frac{\partial u}{\partial x} = 0
	\end{equation}

	By projecting Eq.~(\ref{eq:1dadvection}) onto the the space $\pmb{\Omega}$, and combining Eq.~(\ref{eq:discrete_1D_full})~\&~(\ref{eq:discrete_1D_ad}) result in:
	\begin{equation}
		\frac{\partial \bar{\mathbf{u}}_j}{\partial t} = -J_j^{-1} \Big ( \mathbf{D}\mathbf{f}_j + \big (f^{c}_l  - f_j(x_j) \big )\mathbf{h_l} + \big (f^{c}_r - f_j(x_{j+1})\big)\mathbf{h_r} \Big )\label{eq:FR_full_anal}
	\end{equation}
	The notation used is compatible with that set out initially in~\cite{Vincent2011}, taking $\mathbf{D}_{mn}$ as the contribution from the first derivative on the $m^{th}$ Lagrange basis function to the $n^{th}$ solution point. $\mathbf{h_l}$ and $\mathbf{h_r}$ are taken as $dh_l(\pmb{\xi})/d\xi$ and $dh_r(\pmb{\xi})/d\xi$ respectively, where $\pmb{\xi}$ are the reference coordinates of the solution points. By setting $a=1$ and applying upwinding at the interfaces:
	\begin{align}
		J_j^{-1}\hat{f}^{c}_l &= J_{j-1}^{-1}\hat{u}_{j-1}(1) \label{eq:int_jac1}\\		
		J_j^{-1}\hat{f}^{c}_r &= J_j^{-1}\hat{u}_j(1) \label{eq:int_jac2}
	\end{align}		
	By substituting Eq.~(\ref{eq:int_jac1})~\&~(\ref{eq:int_jac2}) into Eq.~(\ref{eq:FR_full_anal}) and collecting the matrix operators into $\mathbf{C}_0$  and $\mathbf{C}_{-1}$: 
	\begin{align}\label{eq:DiscreteFR2}
		\frac{\partial \bar{\mathbf{u}}_j}{\partial t} &= -J_j^{-1} \mathbf{C}_0 \mathbf{u}_j - J_{j-1}^{-1} \mathbf{C}_{-1}\mathbf{u}_{j-1} \\		
		\mathbf{C}_0 &= \mathbf{D} - \mathbf{h_l}\mathbf{l_l}^T\\		
		\mathbf{C}_{-1} &= \mathbf{h_l}\mathbf{l_r}^T
	\end{align}
	where $\mathbf{l_l}$ and $\mathbf{l_r}$ are again compatible with \cite{Vincent2011} and defined such that $\mathbf{l_l}_i$ is the contribution of the $i^{\mathrm{th}}$ Lagrange basis function evaluated at the left interface. $\mathbf{l_r}$ is similarly defined.
	
	Defining the continuous input as a Bloch wave, and projecting onto the discrete solution domain:
	\begin{align}
		u(x,t) &= v \exp{\big(i(kx - \omega t)\big)} \label{eq:bloch1d2} \\
		\mathbf{u}_j &= \mathbf{v}_j \exp{\big(ik(0.5(\xi + 1)\delta_j + x_j - ct)\big)} \label{eq:bloch1dd2}
	\end{align}
	Inputting this result into Eq.~(\ref{eq:DiscreteFR2}), and setting $\delta_j = x_{j} - x_{j-1}$ gives:
	\begin{equation}\label{eq:FR_Mod_Wavenumber}
		c(k) \mathbf{v} = -\frac{i}{k}\Big (J_j^{-1}\mathbf{C}_0 + J_{j-1}^{-1}\mathbf{C}_{-1}\exp{\big(-ik\delta_{j}\big)} \Big )\mathbf{v}
	\end{equation}
		
	Equation~(\ref{eq:FR_Mod_Wavenumber}) shows that the modified phase velocity $c(k)$, is one of the complex eigenvalues of a matrix describing the spatial transformation performed by the scheme. For an FR scheme with order $p$ there will be $p$ eigenvalues to this problem, of which one is the physical result and the other modes being phase shifted values to give an orthogonal set. The physical interpretation of $c(k)$ is that a wave number's dispersion factor is $\Re (c(k))$ and its dissipation factor is $\Im (c(k))$.
	
	The special case that $\delta_j = $~\emph{const} implies that $J^{-1} = $~\emph{const} which was the case investigated in \cite{Asthana2014,Vincent2011}. However, the more general form of Eq.~(\ref{eq:FR_Mod_Wavenumber}) allows  von Neumann analysis to be performed on stretched meshes. Importantly, Eq.~(\ref{eq:FR_Mod_Wavenumber}) shows that the stencil of cells affecting the dissipation and dispersion of an upwinded FR scheme is just the current cell and its immediate neighbour and hence only the local expansion rate is important for behaviour. This is clearly not the case for finite difference schemes above second order. Repeating the analysis for an FD scheme will give a  basis of comparison and an example of the modified wavenumber for a $4^{\mathrm{th}}$ order central difference scheme is given in Eq.~(\ref{eq:FDCD4_Mod_wavenumber}):
	\begin{align}\label{eq:FDCD4_Mod_wavenumber}
		c(k) = \frac{i}{k} \Big( &b_{-2,j}\exp{\big(-ik(\delta_j + \delta_{j-1})\big)} +                 b_{-1,j}\exp{\big(-ik(\delta_j)\big)} + \nonumber \\
		&b_{2,j}\exp{\big(ik(\delta_{j+2} + \delta_{j+1})\big)} + b_{1,j}\exp{\big(ik(\delta_{j+1})\big)} \Big)
	\end{align}
	where $b_{-2,j}$ is a weighting factor from the derivative of the Lagrange polynomial basis function corresponding to the point $j-2$ evaluated at the point $x_j$, and so on.    
	
	A further implication of Eq.~(\ref{eq:bloch1dd2}) is that Eq.~(\ref{eq:FR_full_anal}) can be rewritten in a form called the update equation:
	\begin{equation}
		\frac{\partial \bar{\mathbf{u}}_j}{\partial t} = \mathbf{Q} \bar{\mathbf{u}}_j
	\end{equation}
	
	And for the case of pure upwinding at the interface it follows that $\mathbf{Q} = -J_j^{-1}\mathbf{C}_0 - J_{j-1}^{-1}\mathbf{C}_{-1}\exp{\big(-ik\delta_{j}\big)}$. Putting the result into this form allows the analytical framework of Asthana~\cite{Asthana2014} to be used, hence:
		\begin{align}
			\bar{\mathbf{u}}_j(t + \tau) &= \mathbf{R}(\mathbf{Q})\bar{\mathbf{u}}_j(t) \label{eq:RK_gen} \\
			\mathbf{R}_{33} &= \mathbf{I} + \frac{(\tau\mathbf{Q})^1}{1!} + \frac{(\tau\mathbf{Q})^2}{2!} + \frac{(\tau\mathbf{Q})^3}{3!} \label{eq:RK_33}
		\end{align}	
		
		Equation~(\ref{eq:RK_33}) is an example definition of $\mathbf{R}$ for a $3^{rd}$-order 3-step low storage Runge-Kutta\cite{Kennedy2000}. The form of Eq.~(\ref{eq:RK_gen}) implies the von Neumann condition of $\mathbf{R}$'s spectral radius, $\rho(\mathbf{R}) \leqslant 1 \:\: \forall \: k \in \mathbb{R}$. 
		
	From the analytical form the main derived quantity to be considered is the Points Per Wavelength (PPW) for dispersion $\mathrm{error} < 1\%$ , defined as:
	\begin{equation}\label{eq:PPW}
		\mathrm{PPW} = \frac{2\pi}{\big\{ \hat{k} | \inf{(|\Re{(c(\hat{k}))}-1|_2,\epsilon)} \big\}}
	\end{equation}	 
	where, $\epsilon$ is the error level and $c(k)$ is the convective velocity from Eq.(\ref{eq:FR_Mod_Wavenumber}).  This definition of PPW is based on the points being the solution points in the case of FR. Consider the case of a  awve whoser wavelength is the length of an element. The normalised Nyquist wavenumber would then be $2\pi/(p+1)$ and hence the use solution points is contained within the definition. PPW is particularly interesting as it can be used to produce minimum point requirements for a given region, if the scale of flow features or the explicit filter width is known. A further derived quantity that will be touched on briefly is the implicit filter kernel, Trefethen~\cite{Trefethen1994}:
	\begin{equation}\label{eq:kernel}
		\hat{G}(\hat{k}) \propto e^{t\Im(c(\hat{k}))}
	\end{equation}		
	Although this is not the main subject of this paper it can be illuminating to briefly look at the implicit filter.
	
\begin{figure}
		\centering
		\includegraphics[width=0.5\linewidth]{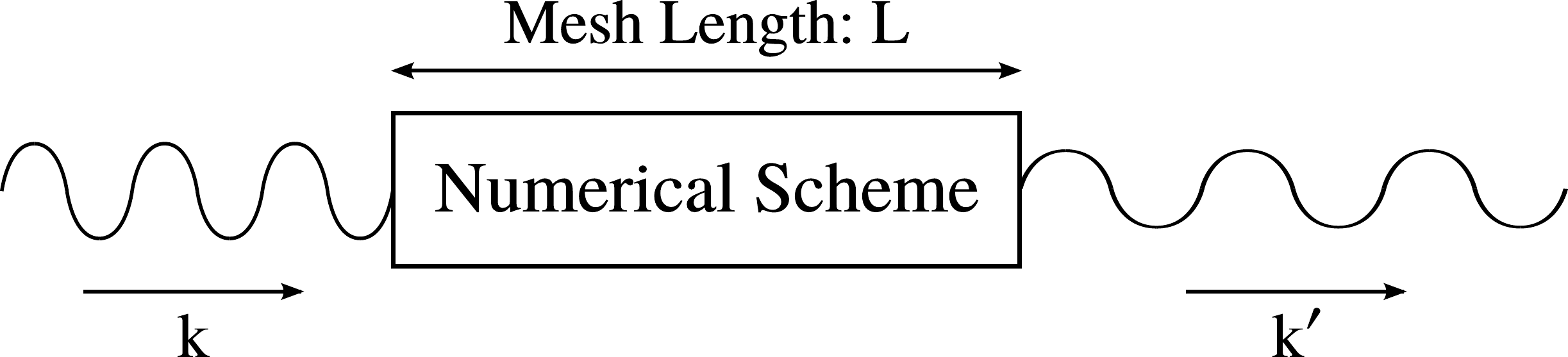}
		\caption{Schematic of numerical scheme with incoming and outgoing wave showing wavenumber transformation over mesh of length $L$.}
		\label{fig:InOutWave}
	\end{figure}

	To validate the analytical methods presented and gain insight into the fully discretised scheme behaviour, a numerical approach is proposed similar to Lele's~\cite{Lele1992} analytical method. The methodology is to apply a scalar wave to a one dimensional domain and allow the wave to be convected downstream. By choosing a low CFL number of $0.01$, the spatial terms dominate the overall error. (CFL numbers of 0.05 and 0.005 were also tested, and the results were found to be largely unchanged, only affected by the discrete Fourier transform as the sampling rate will vary with CFL. Hence a CFL number of 0.01 was used to mitigate this error and give fast test turnover time). Fourier analysis is then performed on the prescribed wave after convection through the grid. The transform of the field $u$ is:
	\begin{equation}\label{eq:fourier_transform}
		\tilde{u}(x) = \sum_{k=-N/2}^{N/2} A_k \exp{\bigg(\frac{2\pi ikx}{L}\bigg)}
	\end{equation}
	The method by Lele~\cite{Lele1992} would have, however, taken the spatial derivative of the wave after convection to give:
	\begin{equation}\label{eq:d_fourier_transform}
		\tilde{u}^{\prime} (x)  = \sum_{k=-N/2}^{N/2} A_k \bigg(\frac{2\pi ik}{L} \bigg ) \exp{\bigg(\frac{2\pi ikx}{L}\bigg)}
	\end{equation}
	By dividing Eq.~(\ref{eq:d_fourier_transform}) by Eq.~(\ref{eq:fourier_transform}) the modified wavenumber can be obtained as in Eq.~(\ref{eq:modWaveEq}). In practice this approach has low throughput, due to need to calculate the derivative and as the full wavenumber space has explored. Furthermore, this method was found to be prone to the introduction of additional sources of error caused by the method used to calculate the derivative. Therefore, a direct comparison is made between the Fourier Transform~(FT) of the prescribed input wave and the FT of the prescribed wave after convection through the grid. This is schematically represented in Fig.~\ref{fig:InOutWave}. 
	\begin{equation}\label{eq:modWaveEq}
		k^{\prime} = \bigg ( \frac{L}{2\pi i} \bigg )\frac{\tilde{u}^{\prime} (x)}{\tilde{u}(x)} 
	\end{equation}
	where $k^{\prime}$ is the modified wave number and equals $\omega/c(k)$. This method is utilised for calculation of the PPW for numerical test cases in accordance with Eq.(\ref{eq:PPW}). 

\section{Analytical Testing}	
	\subsection{Spatial Characteristics}
\begin{figure}
		\centering
		\includegraphics[width=0.45\linewidth,trim= 0mm 0mm 0mm 0mm,clip=true]{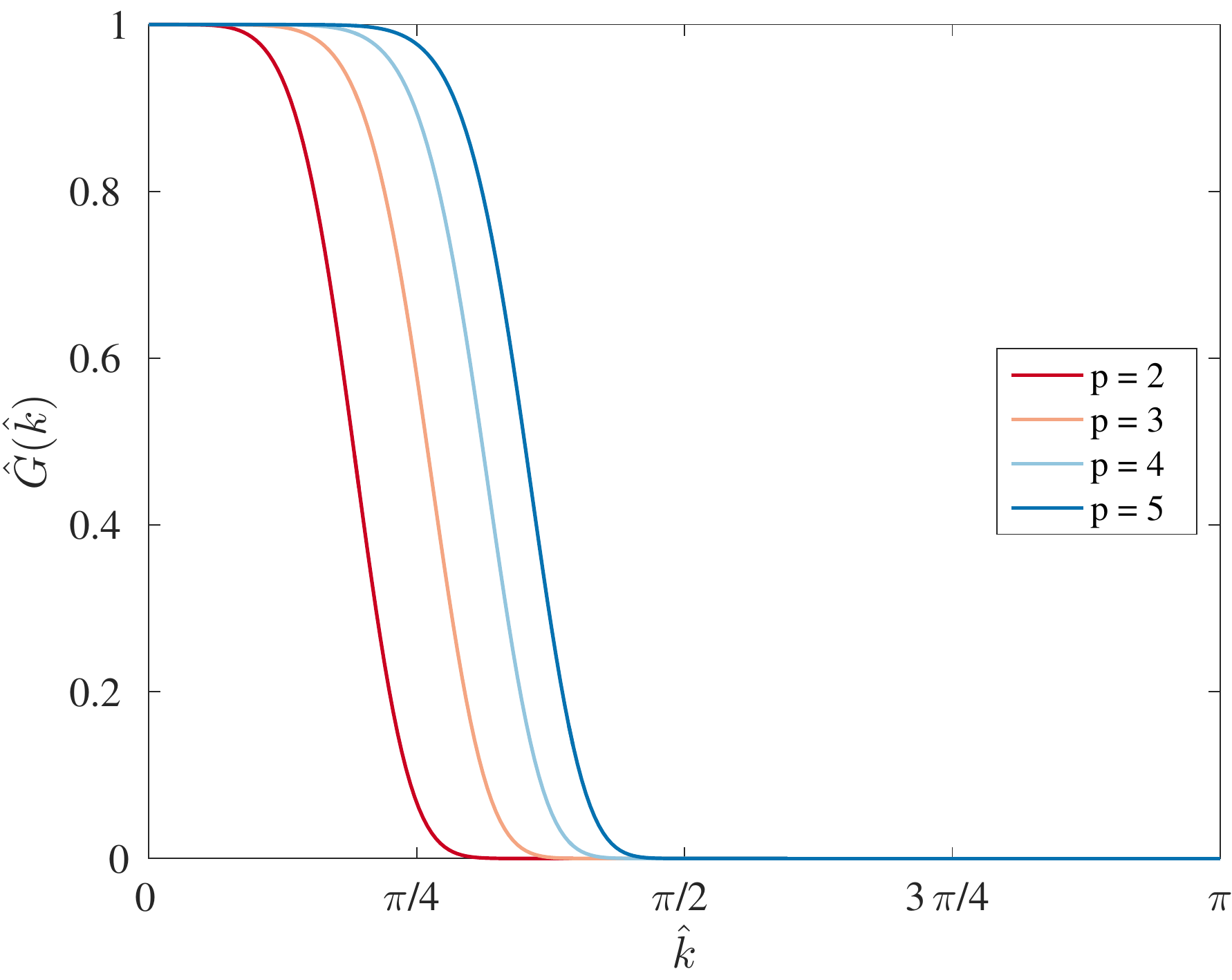}
		\caption{Analytical filter kernel, $\hat{G}(\hat{k})$, against Nyquist normalised wavenumber, $\hat{k}$ for various order of upwinded FR with Huynh, $g_2$, correction functions. The is for $t=100$ in Eq.~(\ref{eq:kernel}).}
		\label{fig:ConvFR180Num}
	\end{figure}

	The analytical method presented in Section~\ref{sec:vn} makes the implicit assumption that there is a linear mapping between solution point placement in the real and physical domain. The analysis did, however, carry through the ability for the relative scaling of adjacent cells to be varied. Hence, this allows for the characteristics of non-uniform grids to be investigated, with the geometric expansion being one such linear transformation in common usage, for example in the meshing of boundary layer. The geometric expansion is defined as:
	 \begin{equation}
	 	x_{j+1} = x_j + \gamma(x_j-x_{j-1})
	 \end{equation}
	 where $\gamma$ is the grid expansion rate. The points that this transformation defines are then used as the flux points for the element and linear interpolation gives the solution points physical domain coordinates, using the computational quadrature as weights to ensure a linear mapping. Proceeding to analytically calculate the modified wave speed from Eq.~(\ref{eq:FR_Mod_Wavenumber}), a preliminary result that can be qualitatively informative is the filter kernel and is shown in Fig.~\ref{fig:ConvFR180Num}. In each case the convolution kernel is normalised independently by the Nyquist wavenumber, $k_{nq}$. This is included as it clearly shows that while there is merit in going to much higher orders ($p>4$), this does result in a case of diminishing returns.
	
\begin{figure}[h]
	 	\centering
		\begin{subfigure}[b]{0.4\linewidth}
			\centering
			\includegraphics[width=\linewidth]{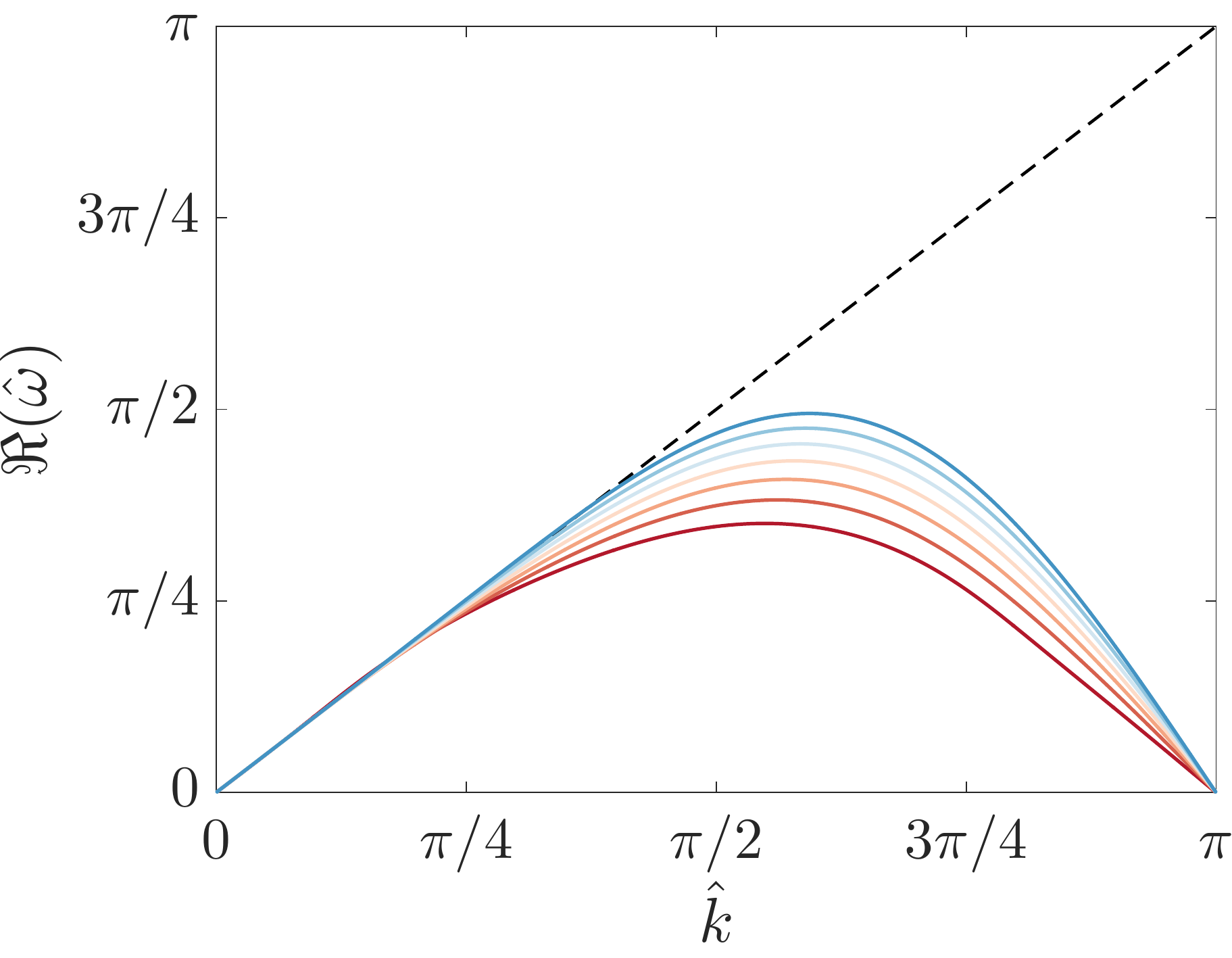}
			\caption{$p=2$}
			\label{fig:stretched_3_R}
		\end{subfigure}
		~
		\centering
		\begin{subfigure}[b]{0.4\linewidth}
			\centering
			\includegraphics[width=\linewidth]{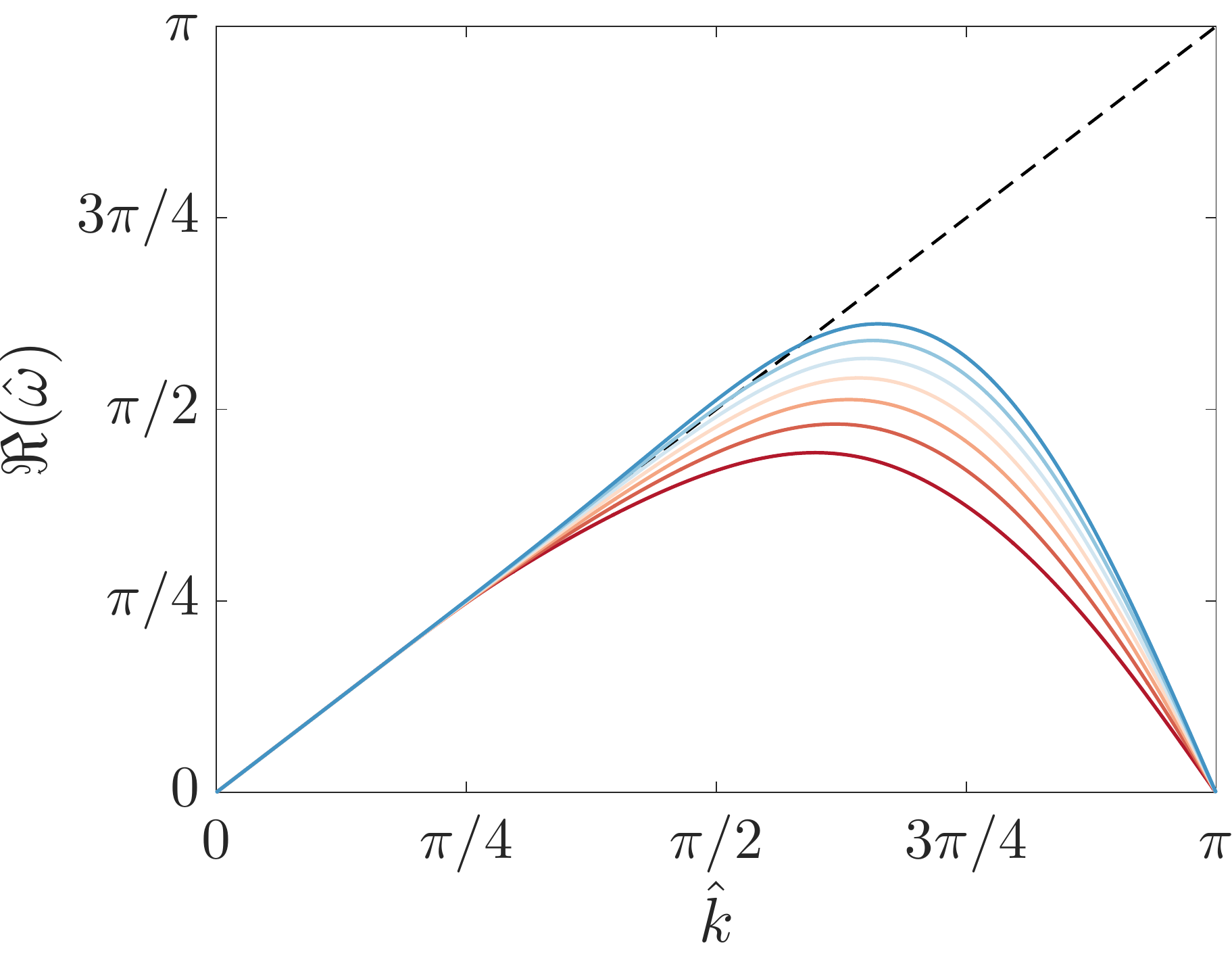}
			\caption{$p=3$}
			\label{fig:stretched_4_R}
		\end{subfigure}
		~
		\begin{subfigure}[b]{0.4\linewidth}
			\centering
			\includegraphics[width=\linewidth]{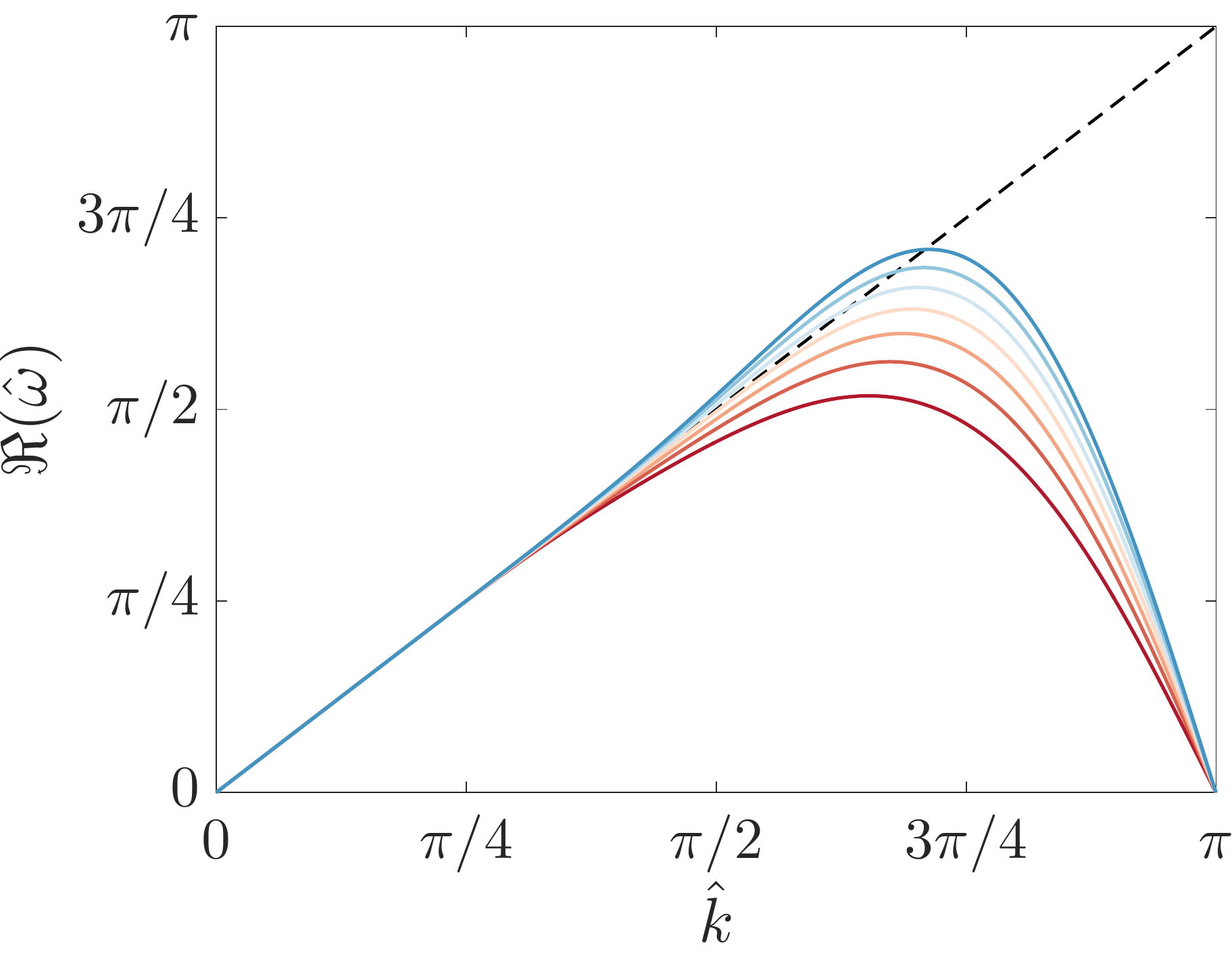}
			\caption{$p=4$}
			\label{fig:stretched_5_R}
		\end{subfigure}
		~
		\begin{subfigure}[b]{0.4\linewidth}
			\centering
			\includegraphics[width=\linewidth]{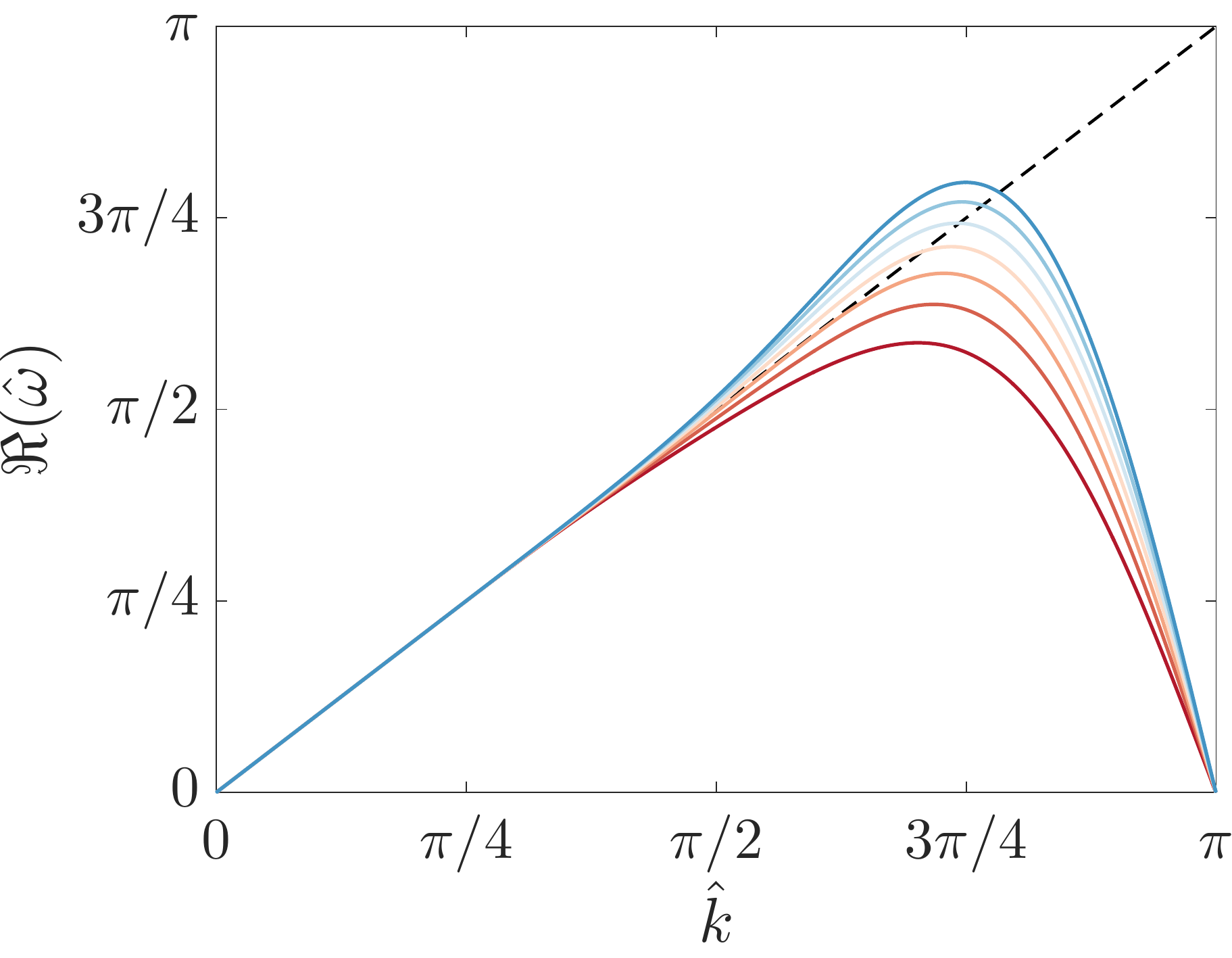}
			\caption{$p=5$}
			\label{fig:stretched_6_R}
		\end{subfigure}
		\caption{Dispersion relations for FR, with Huynh, $g_2$, correction functions, on various geometrically stretched meshes. \\ $--$ spectral performance, {\color{c8_1}\pmb{\pmb{\----}}} $\gamma = 0.4$,{\color{c8_2}\pmb{\pmb{\----}}} $\gamma = 0.6$, {\color{c8_3}\pmb{\pmb{\----}}} $\gamma = 0.8$, {\color{c8_4}\pmb{\pmb{\----}}} $\gamma = 1.0$, {\color{c8_5}\pmb{\pmb{\----}}} $\gamma = 1.2$, {\color{c8_6}\pmb{\pmb{\----}}} $\gamma = 1.4$, {\color{c8_7}\pmb{\pmb{\----}}} $\gamma = 1.6$}
		\label{fig:Stretched_dispersion}
	\end{figure}

\begin{figure}[h]
	 	\centering
		\begin{subfigure}[b]{0.4\linewidth}
			\centering
			\includegraphics[width=\linewidth]{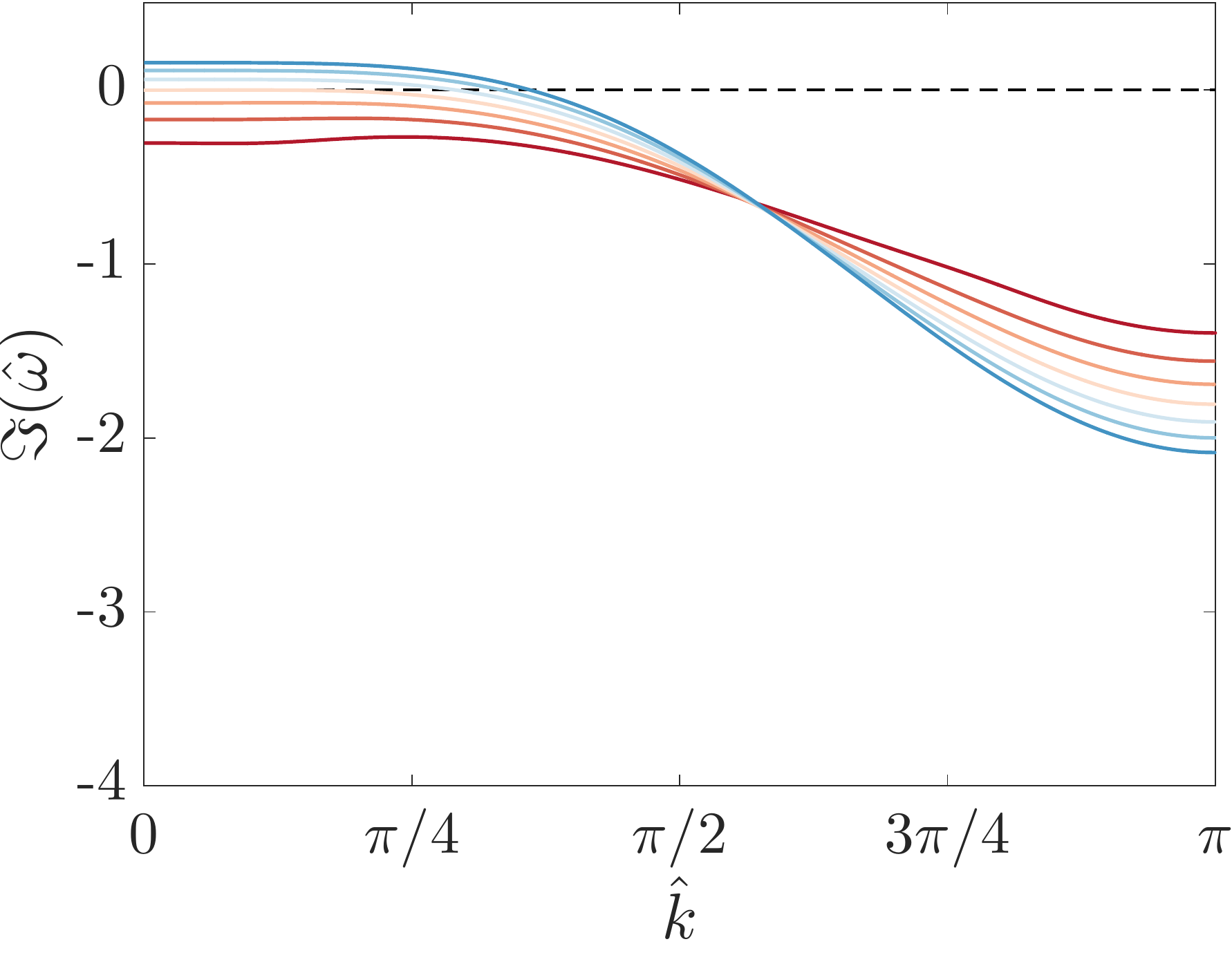}
			\caption{$p=2$}
			\label{fig:stretched_3_I}
		\end{subfigure}
		~
		\centering
		\begin{subfigure}[b]{0.4\linewidth}
			\centering
			\includegraphics[width=\linewidth]{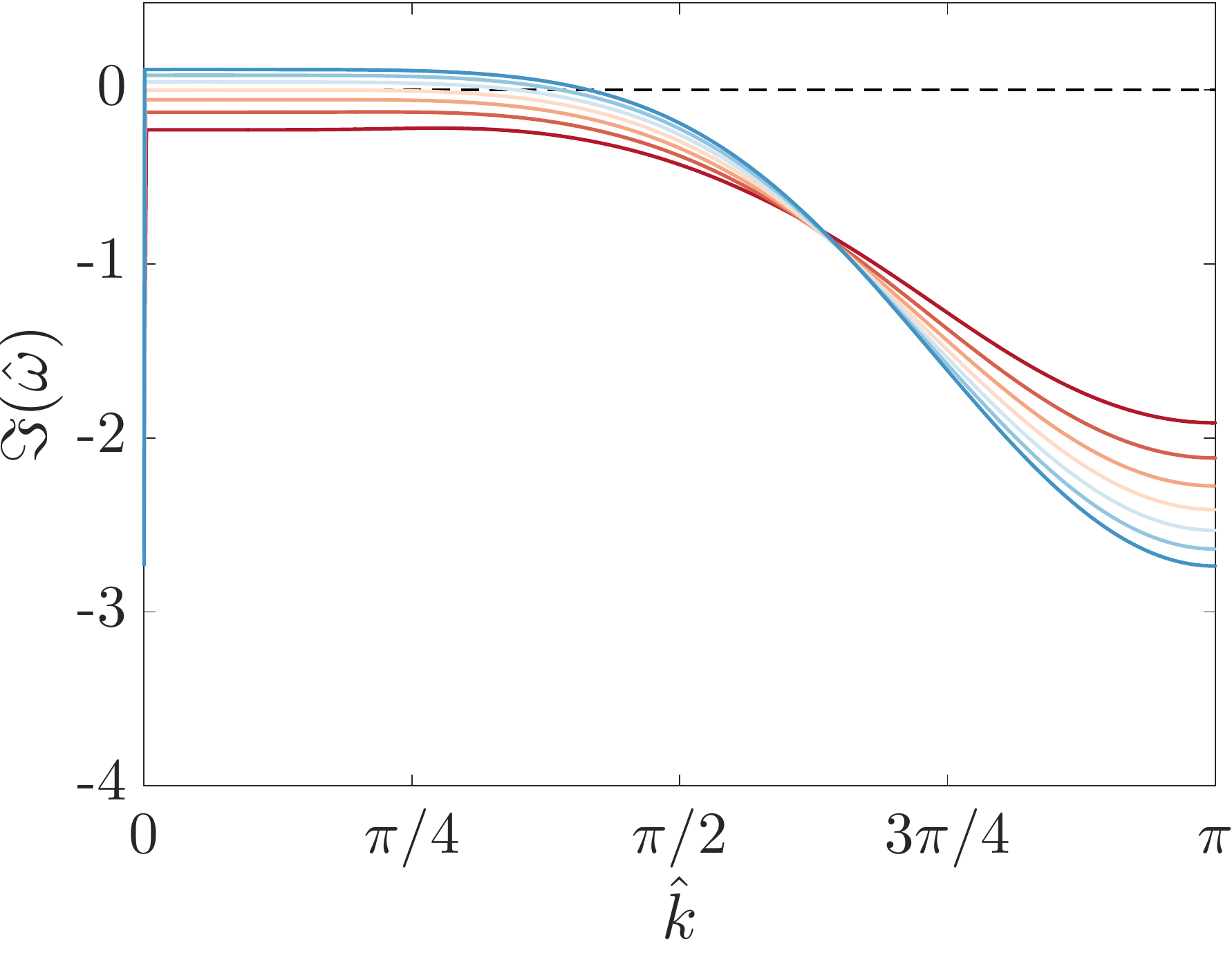}
			\caption{$p=3$}
			\label{fig:stretched_4_I}
		\end{subfigure}
		~
		\begin{subfigure}[b]{0.4\linewidth}
			\centering
			\includegraphics[width=\linewidth]{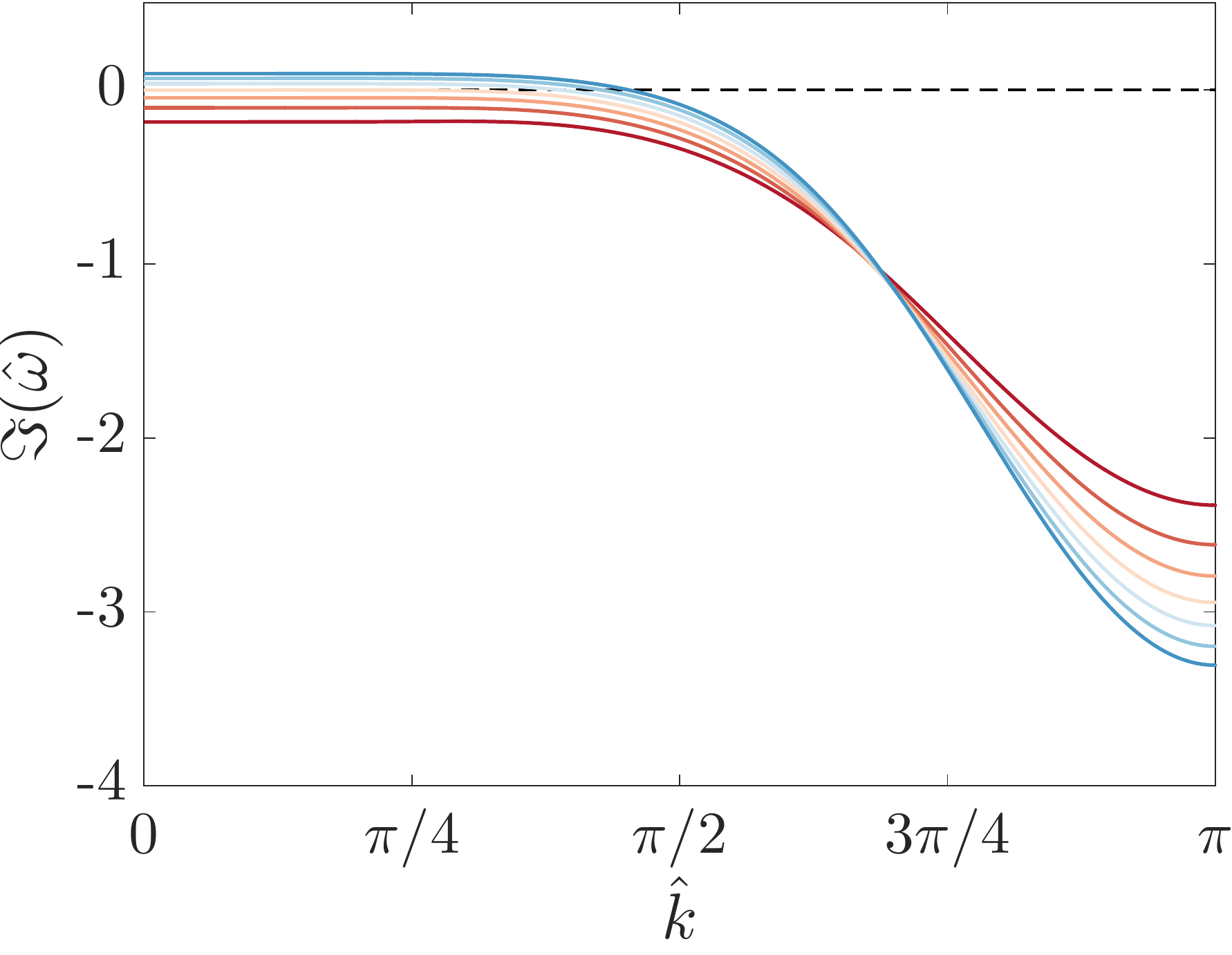}
			\caption{$p=4$}
			\label{fig:stretched_5_I}
		\end{subfigure}
		~
		\begin{subfigure}[b]{0.4\linewidth}
			\centering
			\includegraphics[width=\linewidth]{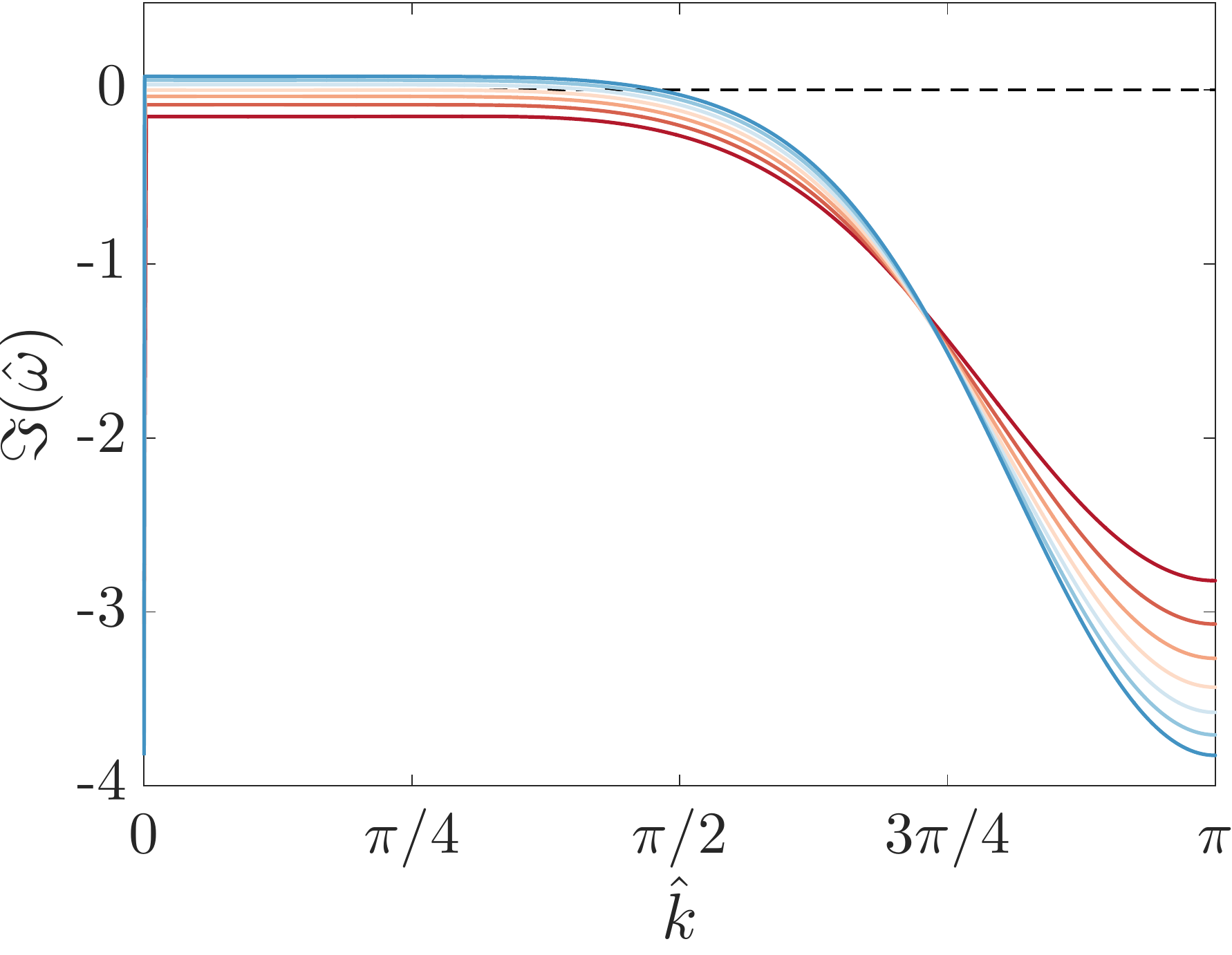}
			\caption{$p=5$}
			\label{fig:stretched_6_I}
		\end{subfigure}
		\caption{Dissipation relations for FR, with Huynh, $g_2$, correction functions, on various geometrically stretched meshes.\\ $--$ spectral performance, {\color{c8_1}\pmb{\pmb{\----}}} $\gamma = 0.4$,{\color{c8_2}\pmb{\pmb{\----}}} $\gamma = 0.6$, {\color{c8_3}\pmb{\pmb{\----}}} $\gamma = 0.8$, {\color{c8_4}\pmb{\pmb{\----}}} $\gamma = 1.0$, {\color{c8_5}\pmb{\pmb{\----}}} $\gamma = 1.2$, {\color{c8_6}\pmb{\pmb{\----}}} $\gamma = 1.4$, {\color{c8_7}\pmb{\pmb{\----}}} $\gamma = 1.6$ }
		\label{fig:Stretched_dissipation}
	\end{figure}

	Figure~\ref{fig:Stretched_dispersion} shows normalised wavenumbers against normalised modified wavenumbers for differing levels of mesh stretching. The results for the uniform grid case can be compared to those found by Huynh~\cite{Huynh2007} and found to be in agreement. What the new results  broadly exhibit is dispersion overshoot for expanding grids, while contracting grids cause dispersion undershoot. For expanding meshes this physically means that for a central band of wavenumbers, where the group velocity ($\partial\omega/\partial k$) is not approximately constant,  the upwind group velocity is higher. This is caused by the change in the Nyquist wavenumber between the smaller upwind and larger downwind elements in an expanding mesh. Resultantly, as the solution advects downwind, this small increase in the group velocity between elements means that a wave will be advected into an element faster than it will exit. This gives rise to anti-dissipation, seen in Fig~\ref{fig:Stretched_dissipation}. The opposite behaviour is exhibited in contractions. At wavenumbers above this central band, Lagrange fitting becomes ineffective at sufficiently projecting the prescribed wave into the functional space, and so the dispersion relation goes zero regardless, and the dissipation becomes high.
	
\begin{figure}
		\centering
		\includegraphics[width=0.5\linewidth,trim= 0mm 0mm 0mm 0mm,clip=true]{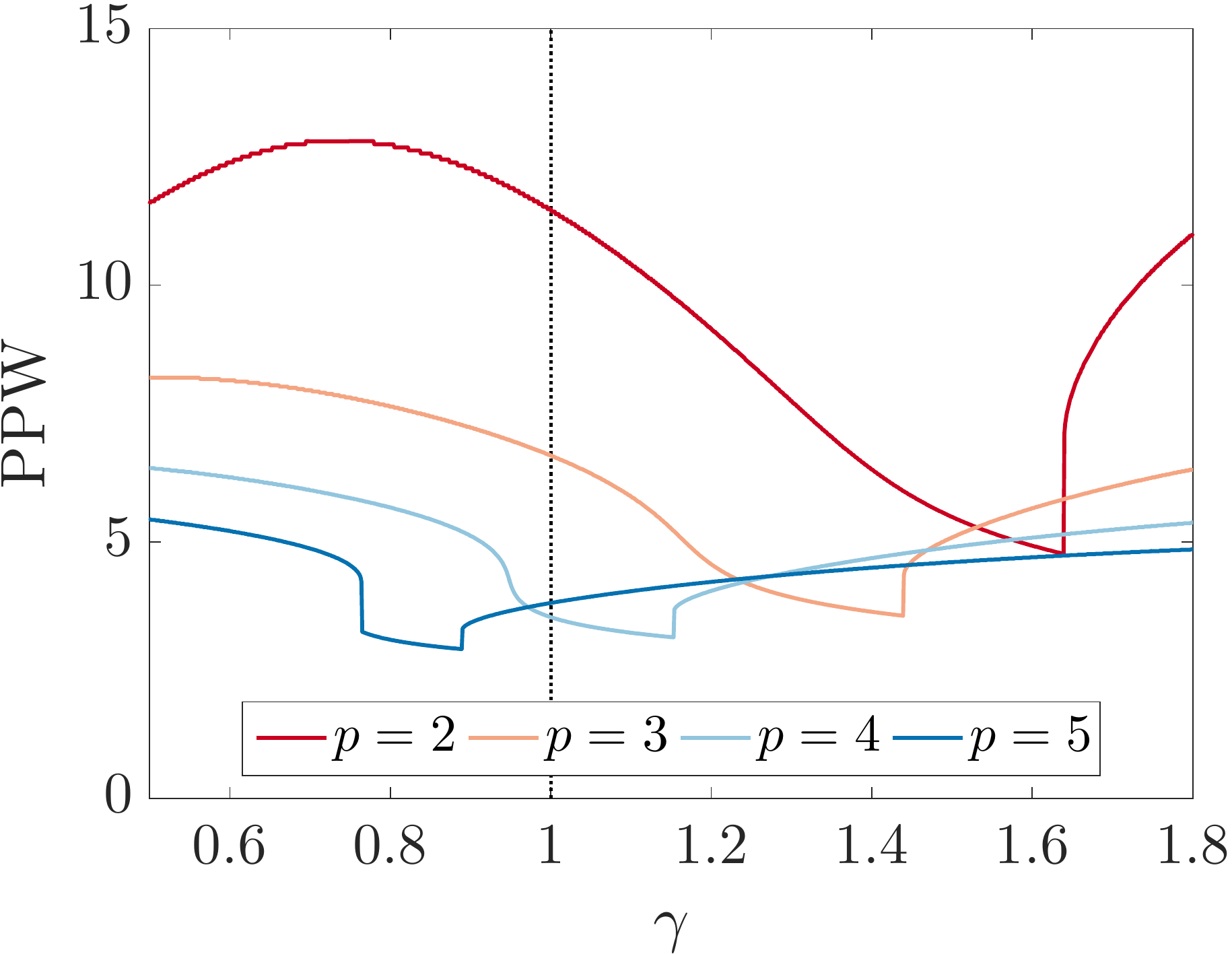}
		\caption{Points Per Wavelength (PPW) for $\mathrm{error} < 1\%$ against expansion rate, $\gamma$, for various spatial orders of upwinded FR with Huynh, $g_2$, correction functions.}
		\label{fig:FR_PPW}
	\end{figure}

	To highlight the practical impact of mesh deformation the PPW resulting in a dispersion $\mathrm{error} < 1\%$, Eq.(\ref{eq:PPW}) with $\epsilon=0.01$, is plotted against expansion factor in Fig.~\ref{fig:FR_PPW}. Over the range of expansion factors there are some clear optimal PPW at varying polynomial orders. When Fig.~\ref{fig:Stretched_dispersion} is considered, it can be seen that the dispersion over- or under- shoot present under  some conditions can be counteracted by mesh warping. Therefore, depending on mesh conditions, it may be beneficial to directionally vary the spatial order, as this may reduce the point requirements locally. For example, for an eddy passing through a complex mesh, fewer points would be needed while passing through a contraction with $p=5$ compared to $p=4$ and vice versa. It is proposed that this could be achieved by maintaining the order of the polynomial interpolation, reflected in the number of solution points, but the order of the correction function could be varied. Clearly this can only be used as a means of dropping the order accuracy, and not as a means of increasing order. From Vincent~\etal~\cite{Vincent2010} it can be seen that this method results in a special case of the energy stability criterion and that the correction functions proposed by Vincent~\etal~\cite{Vincent2010} with reduced order, will still fulfil this criterion. A study was carried  out to this effect, and it was found that using a $p^{\mathrm{th}}$ order, as opposed to $p+1^{\mathrm{th}}$ order, correction function on a $p^{\mathrm{th}}$ order sub-domain results in a stable degradation of the spatial order, pointing to the feasibility of this method.
	
	\subsection{Coupled Spatial-Temporal Characteristics}	
\begin{figure}[h]
		\centering
		\begin{subfigure}[b]{0.324\linewidth}
			\centering
			\includegraphics[width=\linewidth,trim= 1mm 0mm 34mm 0mm,clip=true]{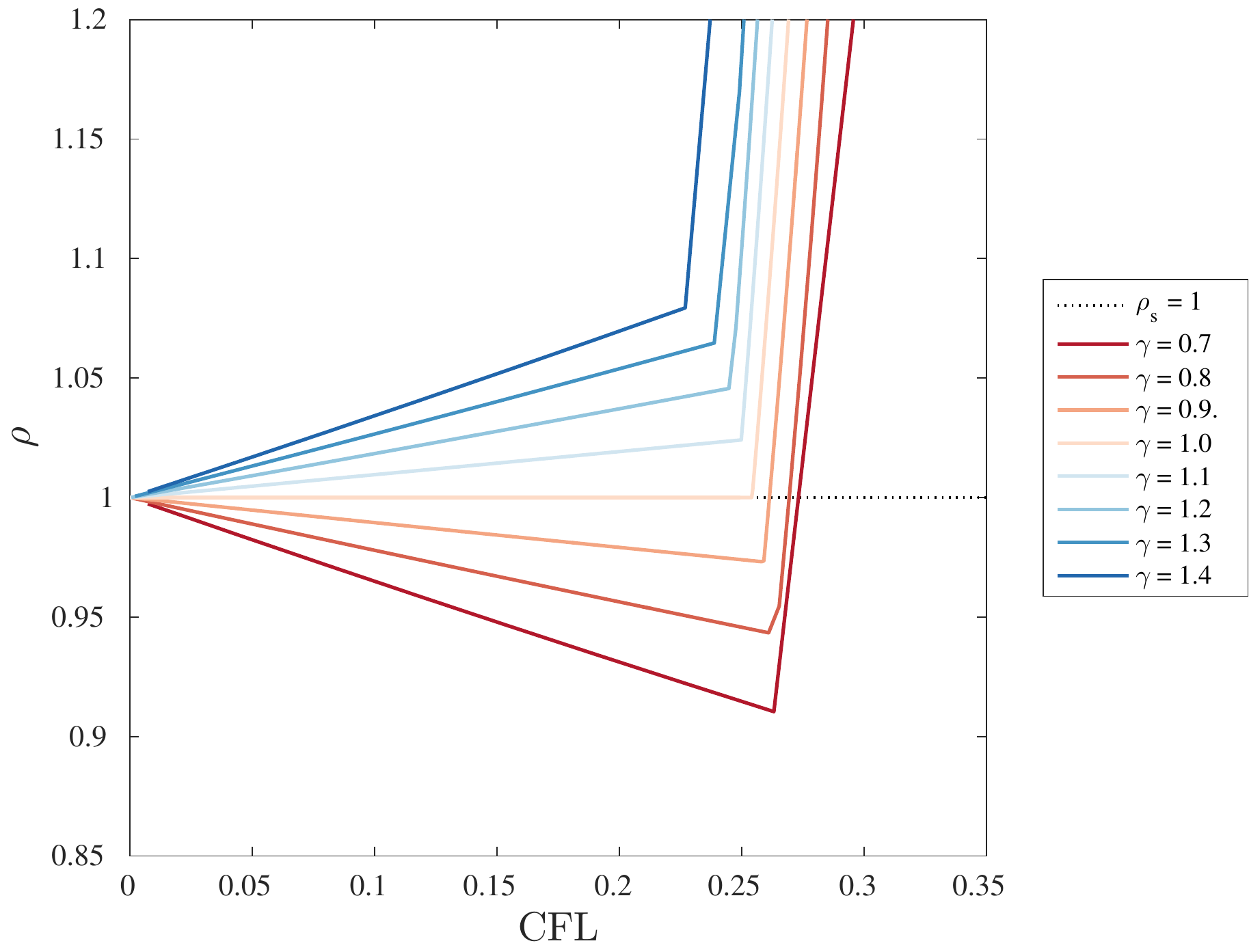}
			\caption{RK33}			
		\end{subfigure}
		~
		\begin{subfigure}[b]{0.4\linewidth}
			\centering
			\includegraphics[width=\linewidth,trim= 1mm 0mm 0mm 0mm,clip=true]{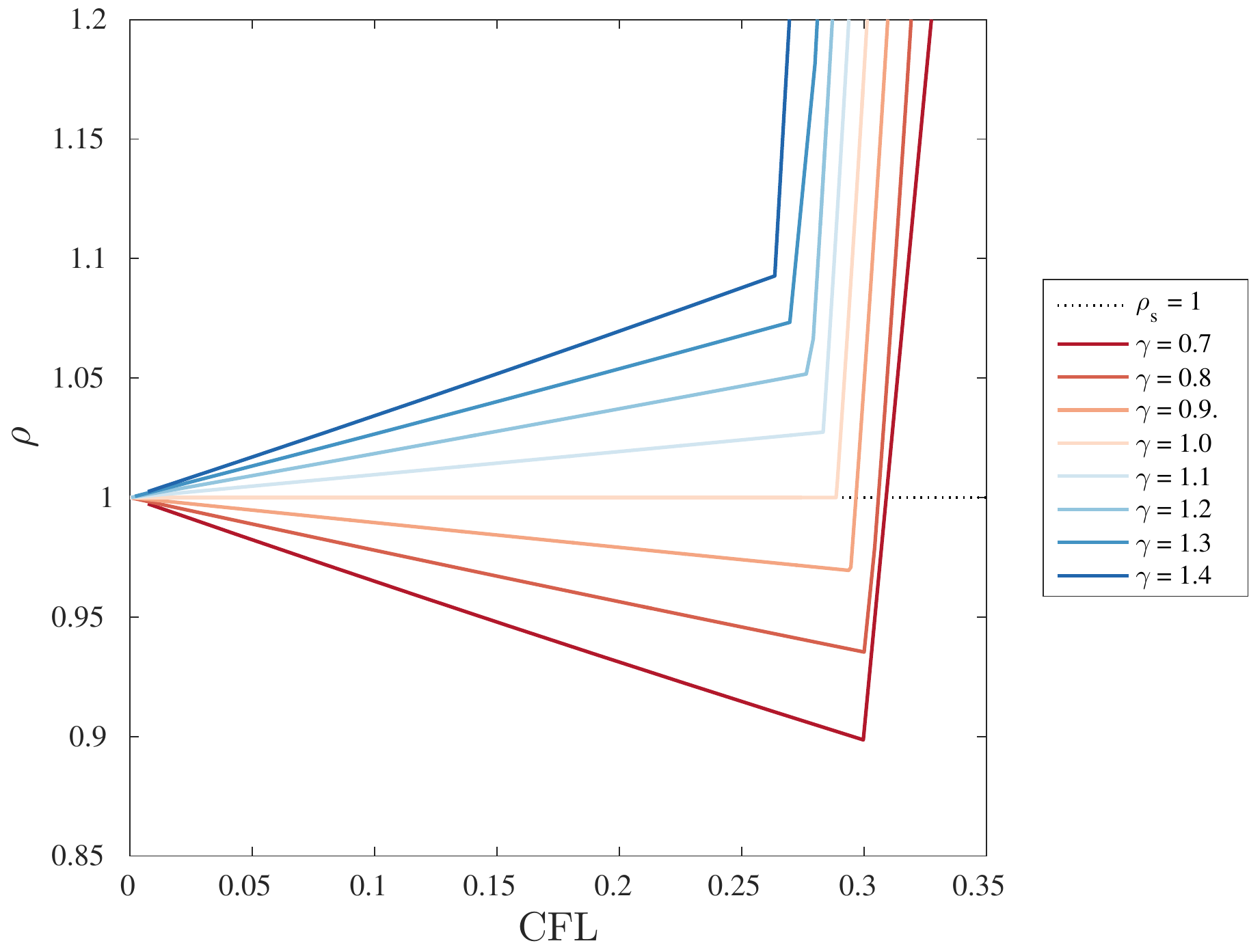}
			\caption{RK44}
		\end{subfigure}		
		\caption{Spectral radius of time scheme specific update matrix, $\mathbf{R}$, for $p=3$ upwinded FR with Huynh, $g_2$, corrections, against CFL number on various grid. }
		\label{fig:RK_spectral_r}
	\end{figure}

	The preceding analysis was predicated on using an analytical solution to the linear advection equation that allowed the time derivative for the semi-discrete linear system to be exactly calculated. This is of course a simplification that cannot be used practically and Eq.(\ref{eq:RK_gen}) is far more representative of a real implementation. The family of temporal integration schemes to be coupled to FR here are low-storage Runge-Kutta~\cite{Kennedy2000}, commonly called RK33, RK44, etc.. The analysis performed is primarily in search of the  maximum stable CFL number, $\tau/\delta_j$ for $a =1$, obtained via varying the time step, $\tau$, and calculation of the spectral radius of the update matrix, $\mathbf{R}$. This is plotted in Fig.~\ref{fig:RK_spectral_r} for various geometric expansion ratios. Initially focusing on contracting grids, the maximum stable CFL number is shown to be higher than in the case of uniform grids. This could have been expected from observation of Fig.~\ref{fig:Stretched_dissipation}, and furthermore it can be reasoned that as a wave is swept from one cell to its smaller upwind neighbour, the ability of the neighbouring element to resolve that wave improves. This is due to the wave's Nyquist normalised wavenumber, $\hat{k}$, decreasing as it is advected through successively smaller elements.
	
\begin{figure}[h]
		\centering
		\begin{subfigure}[b]{0.32\linewidth}
			\centering
			\includegraphics[width=\linewidth]{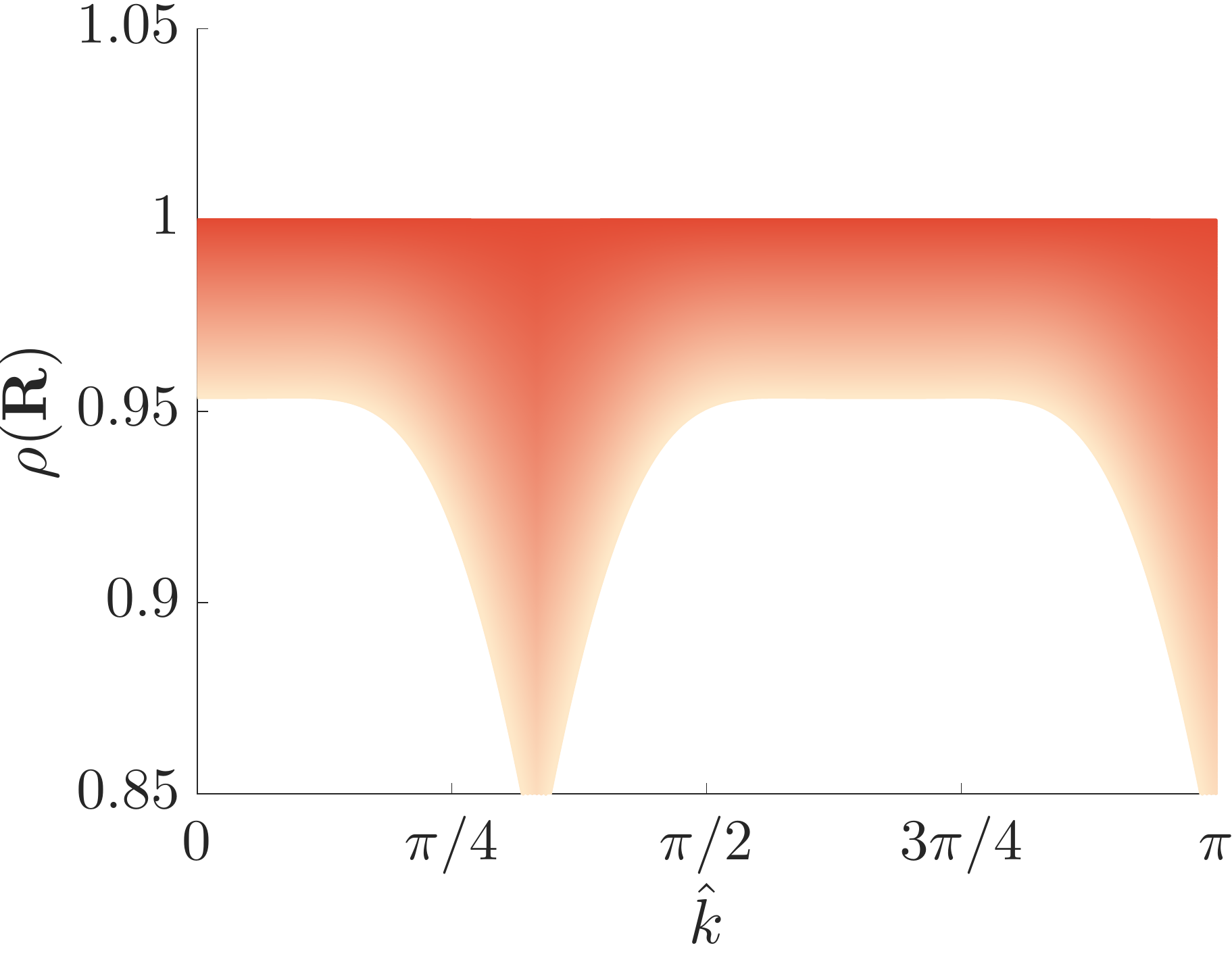}
			\caption{$\gamma = 0.9$, $p=2$}			
			\label{fig:spec_k0.9_44}
		\end{subfigure}
		~
		\begin{subfigure}[b]{0.32\linewidth}
			\centering
			\includegraphics[width=\linewidth]{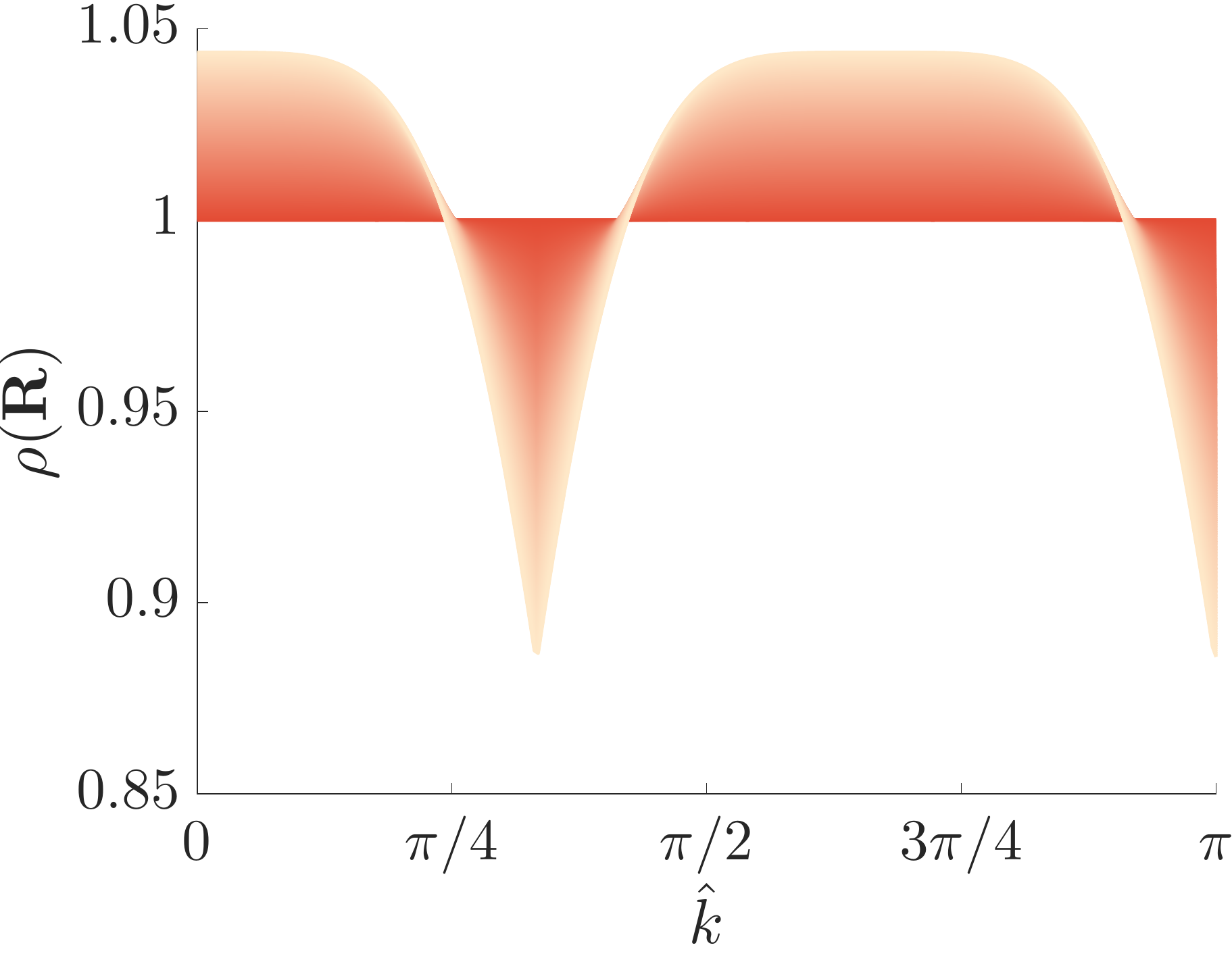}
			\caption{$\gamma = 1.1$, $p=2$}
			\label{fig:spec_k1.1_44}
		\end{subfigure}
		~
		\begin{subfigure}[b]{0.32\linewidth}
			\centering
			\includegraphics[width=\linewidth]{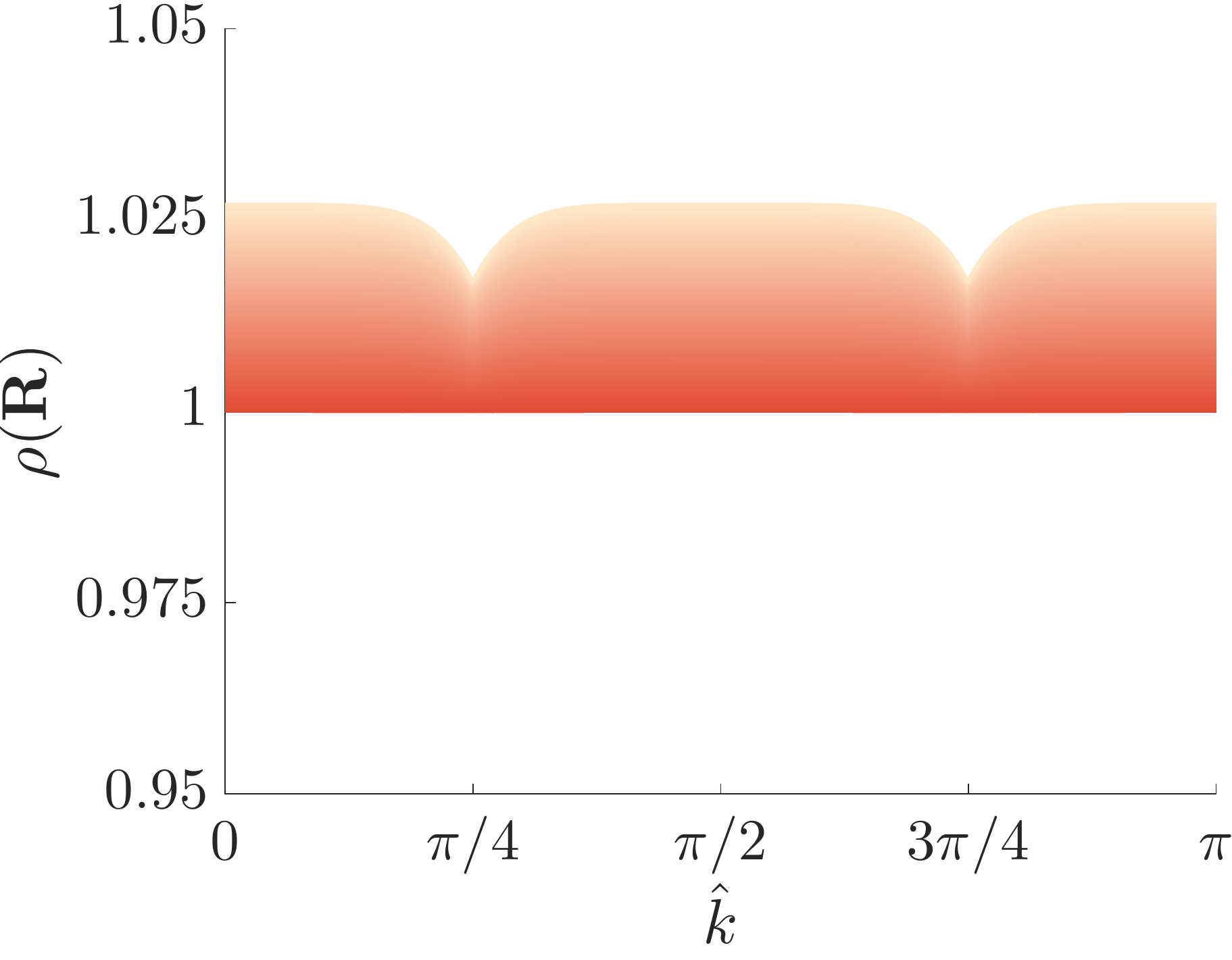}
			\caption{$\gamma = 1.1$, $p=3$}
			\label{fig:spec3_k1.1_44}
		\end{subfigure}			
		\caption{Spectral radius of $\mathbf{R}$ for RK44, upwinded FR with Huynh, $g_2$, correction functions, against Nyquist normalised wavenumber for $\gamma = 0.9$ and $\gamma = 1.1$. with various time steps $\tau$. $\tau$ increasing is shown as a decrease in colour intensity.}
		\label{fig:spec_k}
	\end{figure}

	Focusing on expanding grids, if third order is considered, $p=2$, $\rho(\mathbf{R}) \nless 1 \:\: \forall \:\: \{k \in \mathbb{R} : \gamma > 1\}$, this means that while being strictly unstable, some wavenumbers are in practice stable. This is displayed in Fig.~\ref{fig:spec_k1.1_44} with $\rho(\mathbf{R}(k))$ being both less than and greater than one. The practical implication is that a wave, $k$, fed into the expanding grid can cause an instability if $\rho(\mathbf{R}(k)) > 1$, however, as the wave advects the relative wavenumber increases due to a decreasing $k_{nq}$. Hence a band of  $\rho(\mathbf{R}(k)) < 1$ will be encountered by the wave and the instability will be attenuated. This procedure would be expected to repeat until the wave is beyond the grid resolution.~\footnote{Also of note from Fig.~\ref{fig:spec_k} is that the spectral radius is a periodic function which depends on the element Nyquist wavenumber, rather than the solution point Nyquist wavenumber.}
	
	For $p>2$, different expanding mesh characteristics are seen with $\rho(\mathbf{R}) \geqslant 1 \:\: \forall \:\: \{k \in \mathbb{R} : \gamma > 1\}$, Fig.~\ref{fig:spec3_k1.1_44}. Later numerical tests will show that for $p>2$ instability is also encountered, but this is likewise attenuated as the wave advects. This result can be seen analytically be observing the dissipation relation in Fig.~\ref{fig:Stretched_dissipation} and using the same logic as before. The full impact on stability that this implies will be discussed later alongside numerical findings. However, for simplicity, the scheme stability limit (as shown in Table.~\ref{tab:CFL}) will be taken as the higher point between either the point of the sharp increase in the spectral radius or the point at which the spectral radius increases above 1, (see Fig.~\ref{fig:RK_spectral_r}).		

	\begin{figure}		
		\centering
		\captionof{table}{Analytical CFL limit of FR for various grid expansion factors and temporal integration schemes. Using the Huynh correction function.}
		\begin{tabular}{|c|c||c|c|c|c|c|c|c|}
			\hline
				& & \multicolumn{7}{c|}{CFL} \\ \cline{3-9}
			Time Scheme & Spatial Order & $\gamma = 0.7$ & $0.8$ & $0.9$ & $1.0$ & $1.1$ & $1.2$ & $1.3$ \\ \hline\hline
\multirow{3}{*}{RK33}& 3 & 0.519 & 0.482 & 0.463 & 0.448 & 0.442 & 0.436 & 0.424 \cr\cline{2-9}
				     & 4 & 0.284 & 0.269 & 0.261 & 0.254 & 0.250 & 0.245 & 0.239 \cr\cline{2-9}
				     & 5 & 0.183 & 0.177 & 0.172 & 0.167 & 0.164 & 0.161 & 0.159 \cr\hline\hline
\multirow{3}{*}{RK44}& 3 & 0.592 & 0.547 & 0.531 & 0.513 & 0.505 & 0.495 & 0.507 \cr\cline{2-9}
				     & 4 & 0.318 & 0.307 & 0.297 & 0.288 & 0.282 & 0.278 & 0.270 \cr\cline{2-9}
				     & 5 & 0.218 & 0.199 & 0.194 & 0.189 & 0.186 & 0.182 & 0.179 \cr\hline\hline
\multirow{3}{*}{RK55}& 3 & 0.702 & 0.634 & 0.611 & 0.590 & 0.579 & 0.567 & 0.558 \cr\cline{2-9}
				     & 4 & 0.353 & 0.352 & 0.342 & 0.332 & 0.326 & 0.320 & 0.311 \cr\cline{2-9}
				     & 5 & 0.246 & 0.230 & 0.224 & 0.217 & 0.214 & 0.210 & 0.204 \cr\hline
		\end{tabular}\label{tab:CFL}		
	\end{figure}

\section{Numerical Results}
\subsection{Grid Stretching for Linear Advection}
	The analytical procedures set out up to this point have been semi-discrete and idealised. For CFD practitioners the comparative performance and implemented performance of FR is highly important. From this the process of mesh generation can be informed as well as greater understanding of the expected results gained. To this end, numerical tests are performed for wavenumbers $0\leqslant k \leqslant k_{nq}$, where $k_{nq}$ is the Nyquist wavenumber for a uniform mesh of unit length. For the purposes of comparison, $k$ and $k^{\prime}$ (the modified wavenumber) are normalized by the mesh averaged Nyquist wavenumber. This gives $0 \leqslant \hat{k}= \pi k/k_{nq} \leqslant \pi$ and $0 \leqslant \Re(\hat{k}^{\prime}= k^{\prime}/k_{nq}) \leqslant \pi$. 
	
	Finite Difference (FD) schemes are used to provide a comparison akin to high quality industrial codes. At higher orders, central difference schemes begin to become unstable as, for unstretched grids, second order and greater central difference schemes offer no dissipation. Therefore, the only sources of numerical error are from dispersion and temporal integration. Because of this, and as increasing order leads to better dispersion performance, the cell Reynolds number increases. The result is that FD schemes at intermediate wave numbers become unable to damp out disturbances. (It must be noted that due to the very low CFL number, the temporal scheme introduces negligible numerical error). This is combated by adding smoothing to the solution, a common practice in industry, and is here accomplished by adding a small amount ($0.5-2\%$) of Lax-Friedrichs differencing.
	
\begin{figure}[h]
		\centering
		\begin{subfigure}[b]{0.49\linewidth}
			\includegraphics[width=\linewidth,trim= 32mm 92mm 35mm 98mm,clip=true]{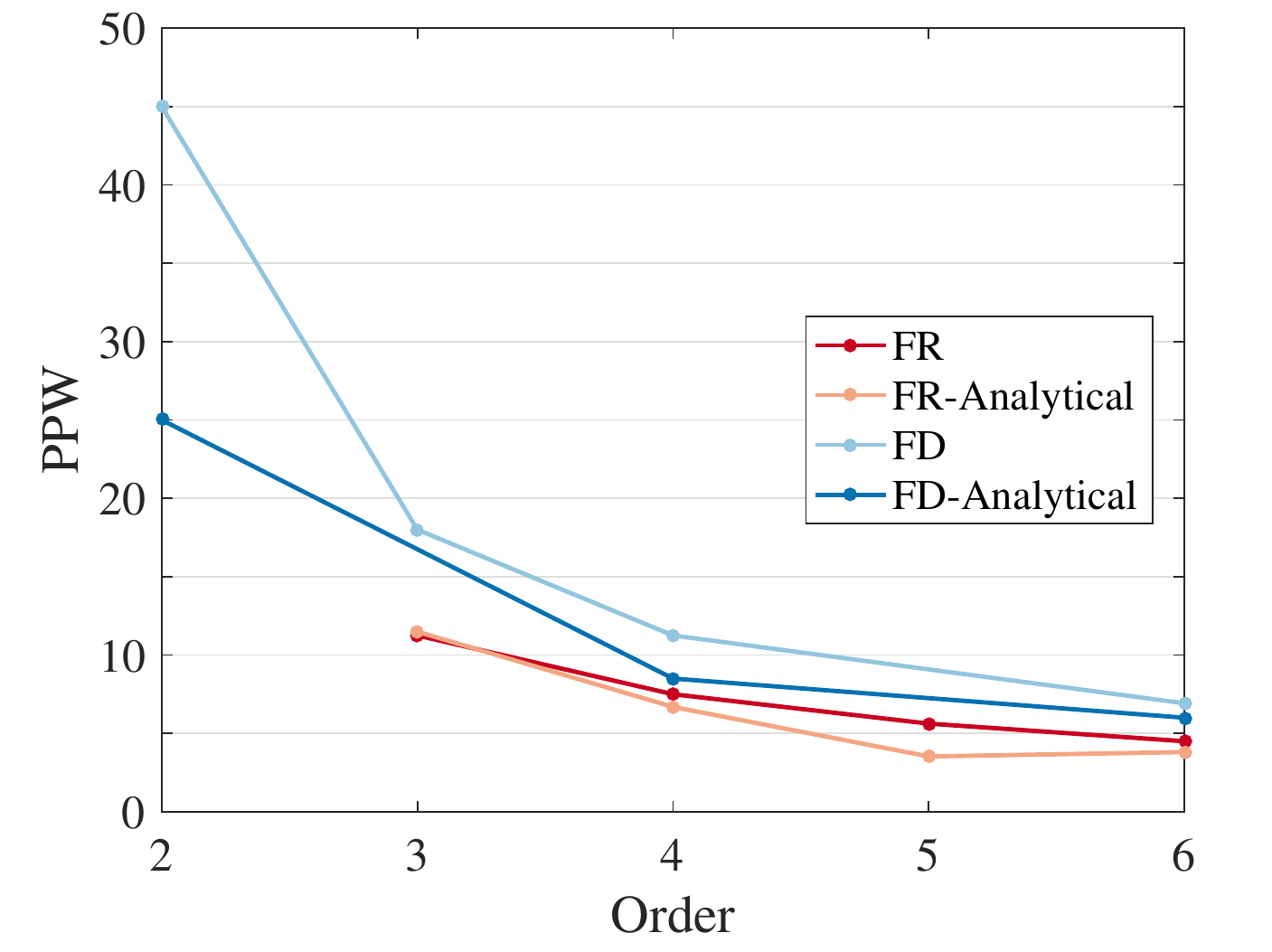}
			\caption{Points per wavelength (PPW) for $1\%$ error against spatial order of FD and FR schemes. FD schemes are central differencing except order 3, which is upwinding. $\mathrm{DoF} = \mathrm{const.} = 180$ and $\tau = 1\times10^{-4}\:s$. The analytical results shown are found from Eq.~(\ref{eq:FR_Mod_Wavenumber}) and similar analysis.}
			\label{fig:FRFD_PPW}
		\end{subfigure}
		~	
		\begin{subfigure}[b]{0.47\linewidth}
			\centering
			\includegraphics[width=\linewidth,trim= 32mm 93mm 42mm 98mm,clip=true]{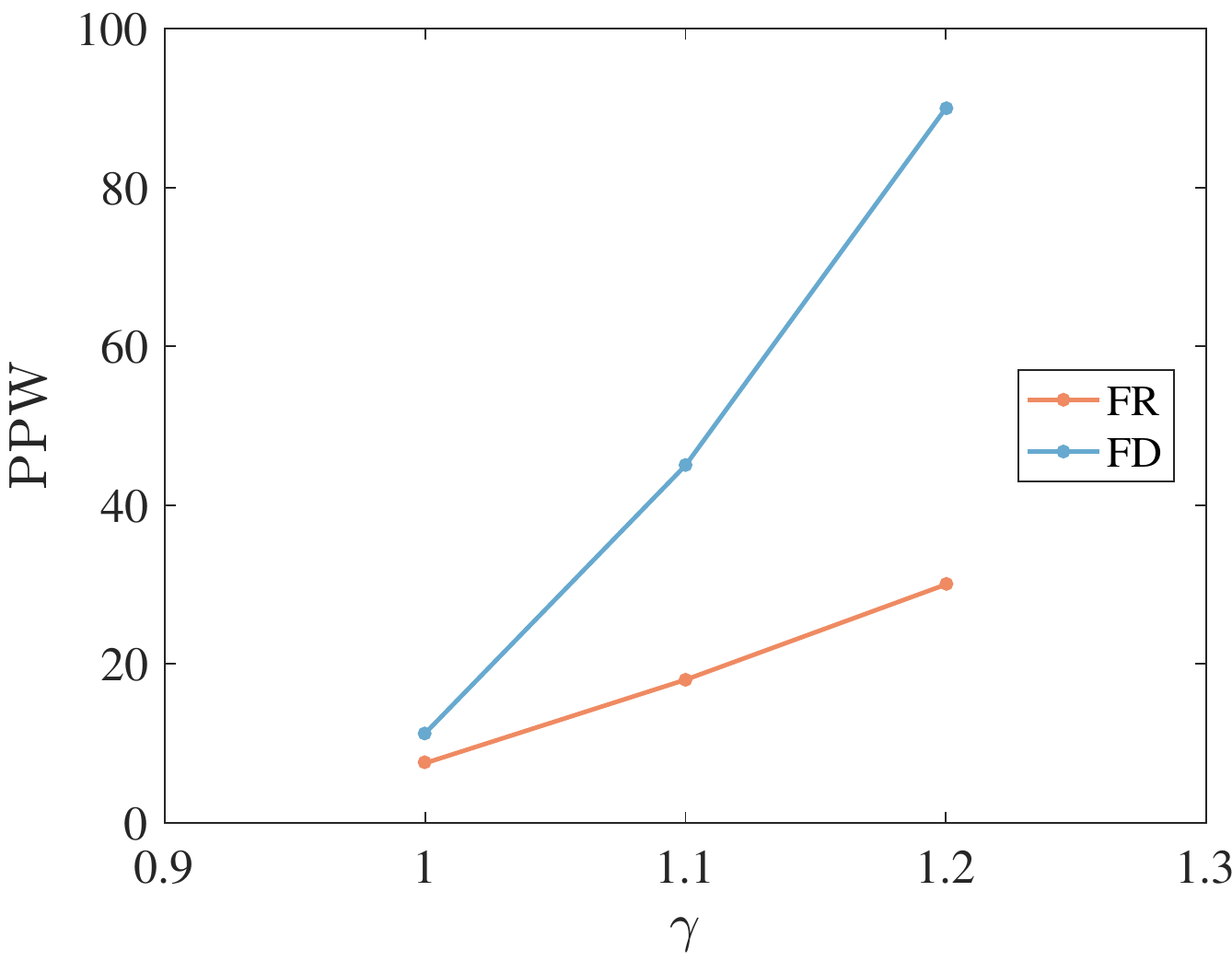}
			\caption{PPW against mesh stretching factor for $4^{\mathrm{th}}$ order FD and FR schemes. $\mathrm{DoF} = 180$  and $\tau =1\times10^{-5}\:s$. \\ \quad \\ \quad \\ \quad \\ \quad }
			\label{fig:PPW_Comp}
		\end{subfigure}
		\caption{}
	\end{figure}

	Two key points are highlighted by Fig.~\ref{fig:FRFD_PPW}. The first of these is that FR requires fewer PPW than FD schemes at equivalent order. Importantly, this means that coarser meshes can be utilised by FR for similar wave resolving chFaracteristics to FD schemes at the same order. It is believed that the increased accuracy of FR originates from the polynomial reconstruction in a reference sub-domain, hence the propagation of information in FR is largely controlled by the correction function, which can lead to superior performance. Whereas, FD methods use a stencil, for which information can freely propagate through, hence a less coherent solution is produced due to each point effectively producing its own polynomial fit of the solution.
	
	Secondly, Fig.~\ref{fig:FRFD_PPW} shows a discrepancy between the theoretical and analytical results for both FR and FD. The origin of the error in both of these schemes is the numerical diffusion. For FD, this is due to the scheme's subtle instability, meaning that for useful implementation, some diffusion must be added to ensure the survival of the solution.  FR is also affected by numerical diffusion, but this is caused by its own intrinsic dissipation, apparent in Fig.~\ref{fig:Stretched_dissipation}. When numerical tests are run, the dispersion and dissipation are inseparable and although a dispersion relationship of sorts can be found, it is impacted by the dissipation of the scheme damping out higher wavenumbers. This is the reality of any application of a scheme and so it can be informative to run both analytical and numerical tests as the numerically derived PPW shown in Fig.~\ref{fig:FRFD_PPW} are those that an end user will experience. The same effect can be seen in Fig.~\ref{fig:PPW_Comp}: as the mesh becomes stretched the PPW rises more quickly than the analytical results would predict. This is caused by the onset of dissipation at lower wavenumbers for deformed meshes, but, importantly, a wave passed through a multi-element mesh will have the transfer function applied multiple times, thereby causing greater attenuation.
	
	Results of great significance displayed both numerically and analytically show that FR has the ability to resolve waves better than FD schemes and that FR is more numerically robust when applied to geometrically stretched meshes, with FR requiring $33\%$ of the mesh points compared to FD in 1D for severely stretched meshes ($\gamma = 1.2$ at fourth order). Moving to two or three dimensions, this result, in the most extreme examples, can be $11\%$ or $4\%$ respectively. The increased ability of FR to handle stretched meshes is again because of the localised fitting within sub-domains, and here the linear transformation caused by the stretching of the elements is exactly captured in the Jacobian. For this case the impact of adjacent cell stretching is felt only through convection of the solution through one interface. However, for a fully compressible Euler or Navier-Stokes implementation the effect of adjacent cells could be increased as Riemann solvers at all interfaces would be necessary and will give rise to more inter-cell communication. 
	
\begin{figure}[h]
		\centering
		\begin{subfigure}[b]{0.44\linewidth}
			\centering
			\includegraphics[width=\linewidth,trim= 0mm 0mm 0mm 0mm,clip=true]{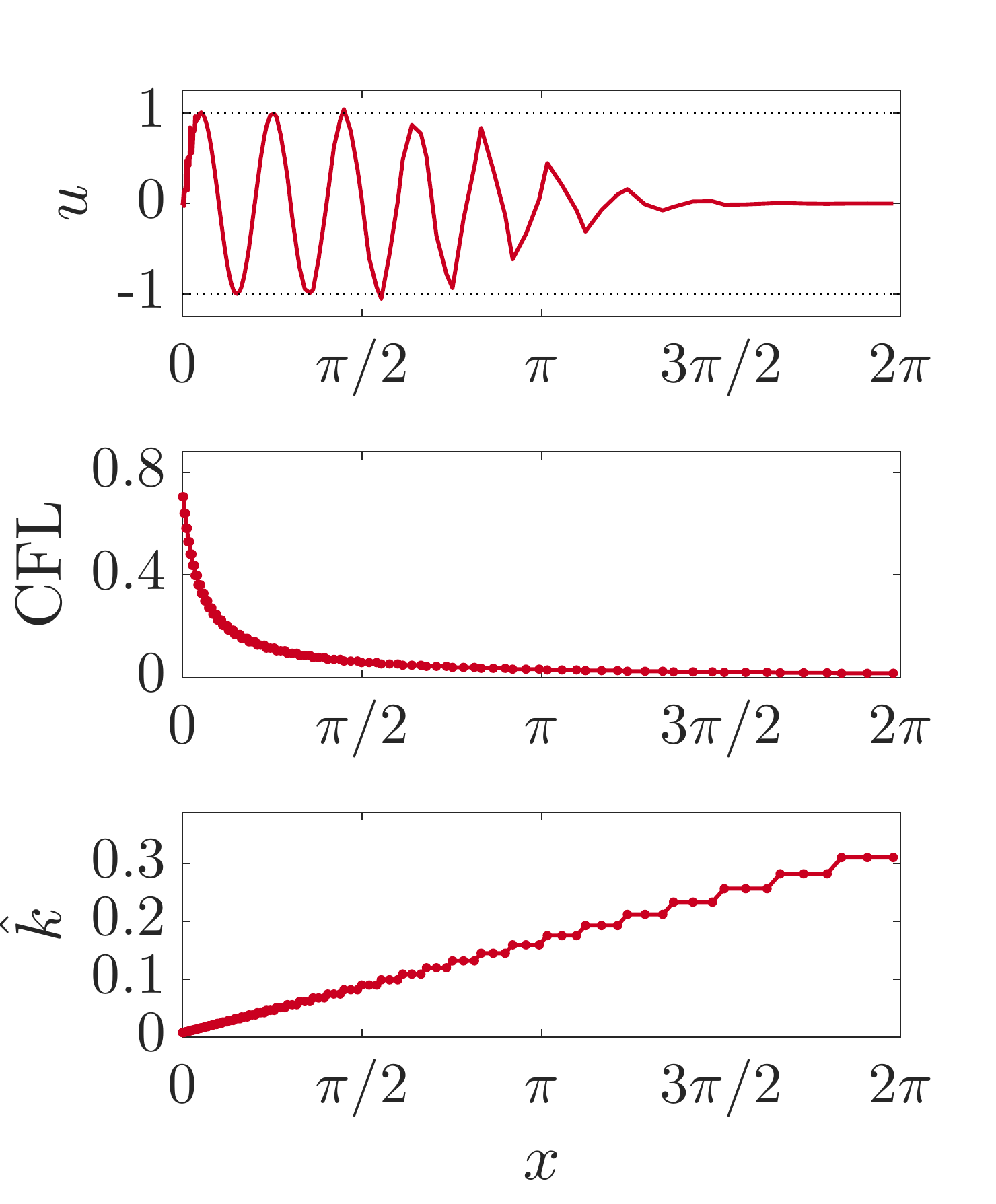}
			\caption{$p=2$}
			\label{fig:CFL_FR3_RK44_11}
		\end{subfigure}
		~
		\begin{subfigure}[b]{0.44\linewidth}
			\centering
			\includegraphics[width=\linewidth,trim= 0mm 0mm 0mm 0mm,clip=true]{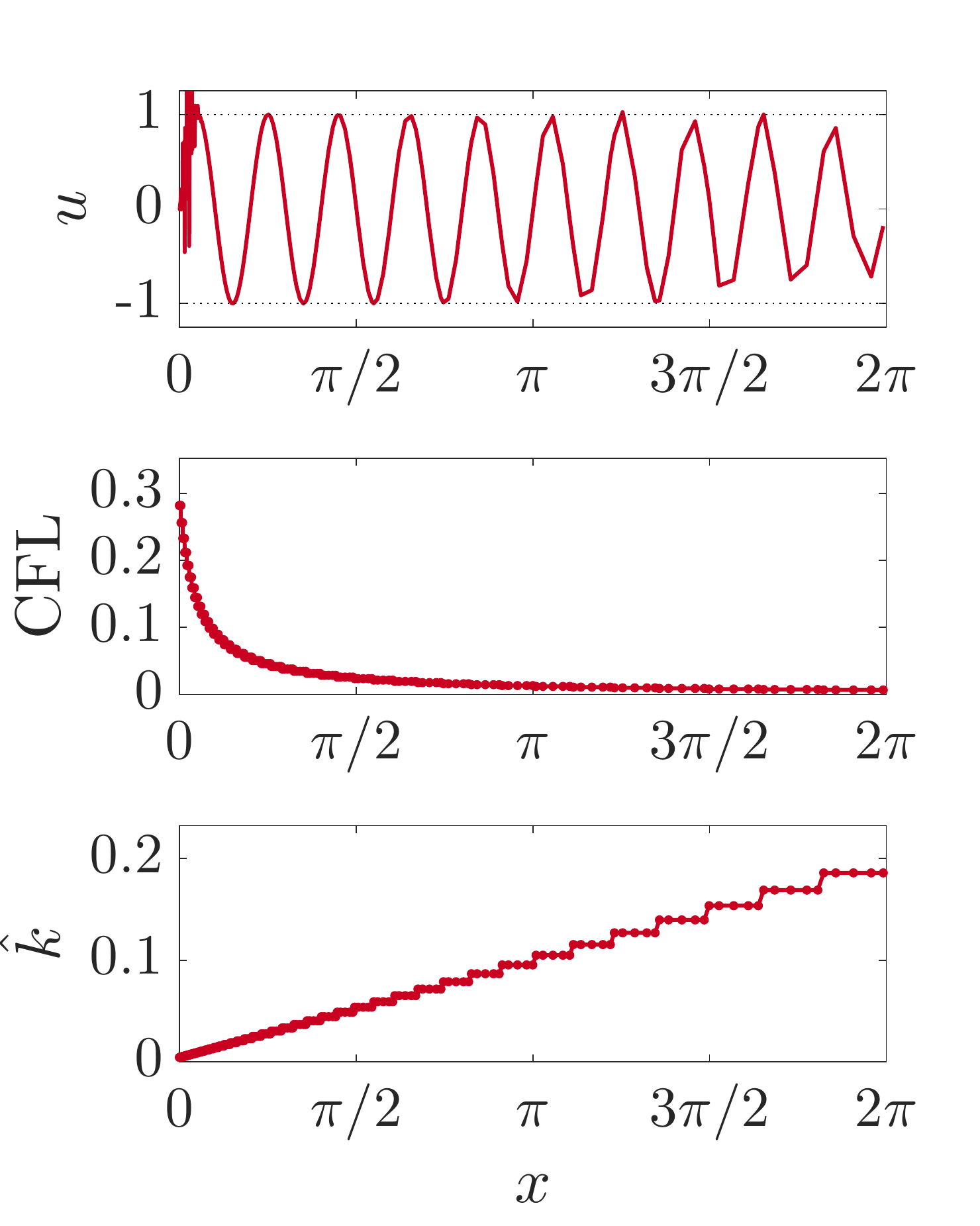}
			\caption{$p=4$}
			\label{fig:CFL_FR5_RK44_11}
		\end{subfigure}		
		\caption{For 3rd order and 5th order FR on an expanding grid with $\gamma = 1.1$ using RK44 time integration a spatial slice is shown. With the convected parameter $u$, the CFL number and the Nyquist normalised ($\hat{k} = \pi k/k_{nq}$).}
		\label{fig:coupled_numerical}
	\end{figure}

	To further understand the stability of the full numerical scheme it is necessary to consider the spatio-temporal coupling. In the previous section the effect of this coupling was considered and it was said that for $p\geqslant3$ on expanding grids, the behaviour is slightly different, \emph{i.e.} $\rho(\mathbf{R}) \geqslant 1 \: \forall \:k$. The implication this has for the stability of higher order expanding grids is not clear from the spectral radius, however, because as the wave moves through the expanding grid $\hat{k}$ will increase, so the scheme dissipation will add a stabilising effect. To show this, a similar numerical method is used, however now taking a spatial slice for various orders, wavenumbers, and grids. Figure~\ref{fig:coupled_numerical} shows two such slices. Initially, the fed wave shows some instability but advection through the grid means dissipation from the spatial scheme will begin to cancel some of the negative dissipation of the time scheme. The result is that after an initial band of instability both orders show recovery of the solution before beginning to decay on the sparser cells. 	
	
	This also illustrated a limitation of the analytical approach adopted here - that taking a solution of the form of Eq.~(\ref{eq:bloch1d2}) means that the solution is static, \emph{i.e.} evolutions of the solution from far upstream are not permitted. Again, this emphasises the importance of running numerical tests alongside analytical ones. 
	
\subsection{Grid Warping for the Euler Equations}
\begin{figure}
		\centering
			\includegraphics[width=0.30\linewidth,trim= 46mm 205mm 90mm 34mm,clip=true]{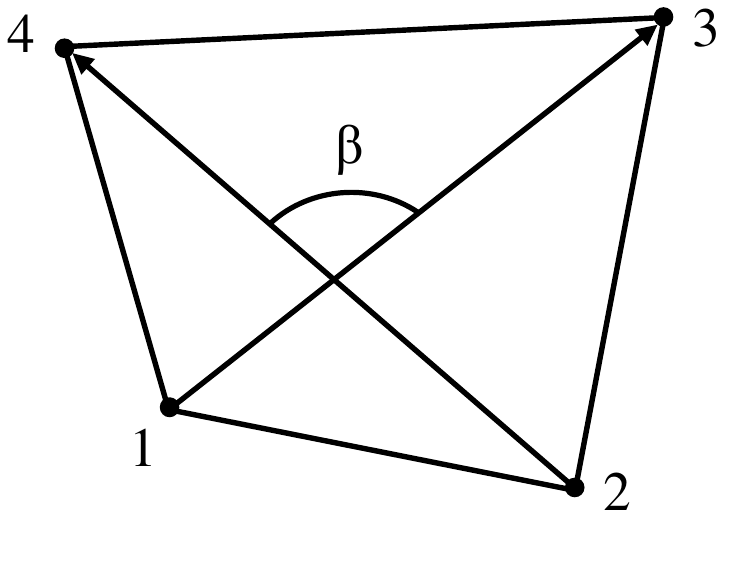}
		\captionof{figure}{Cross diagonal angle definition. $\alpha = \beta - 90^{\circ}$.}
		\label{fig:SkewAngle}
	\end{figure}

\begin{figure}[h]
		\centering
		\begin{subfigure}[b]{0.2\linewidth}
			\centering
			\includegraphics[width=\linewidth,trim= 50mm 100mm 50mm 100mm,clip=true]{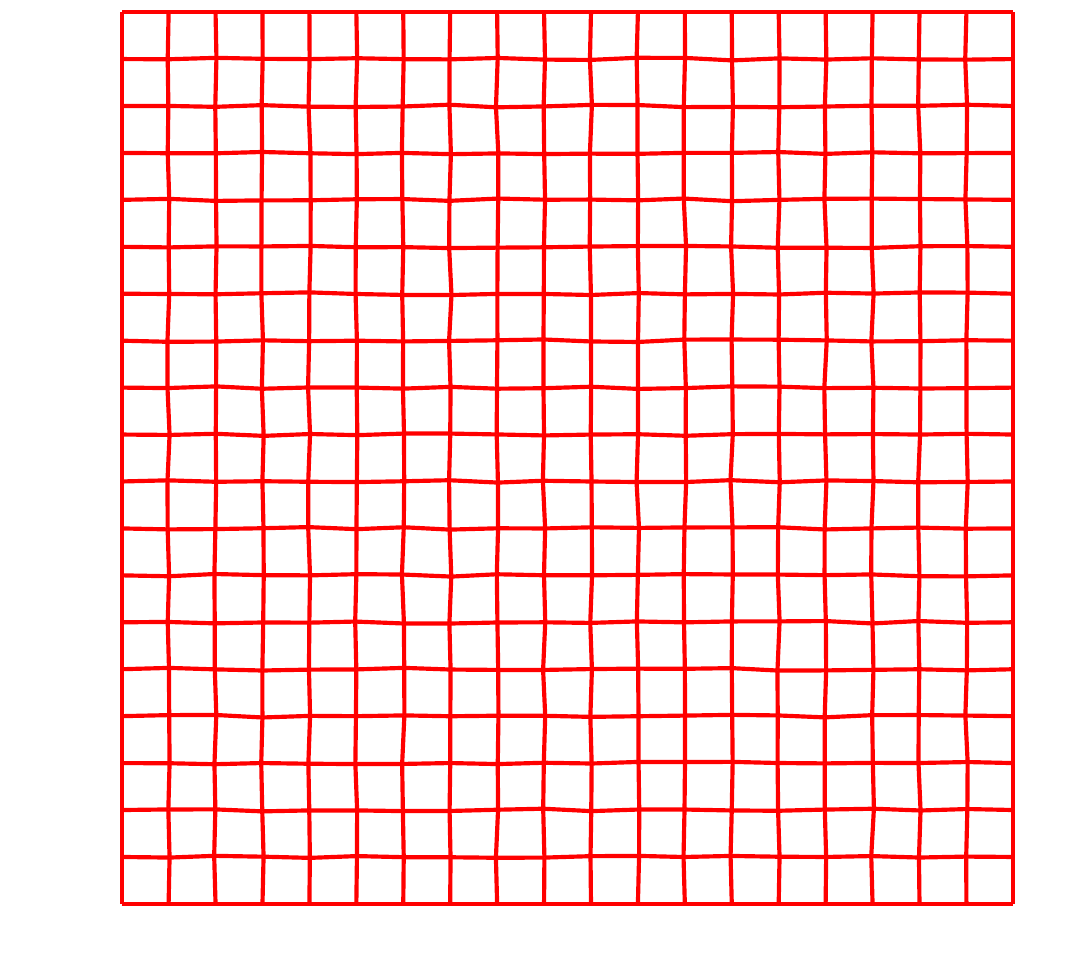}
			\caption{Average skew angle: $\alpha=\:1.0^{\circ}$}
			\label{fig:Jitter1}
		\end{subfigure}
		~
		\begin{subfigure}[b]{0.2\linewidth}
			\centering
			\includegraphics[width=\linewidth,trim= 50mm 100mm 50mm 100mm,clip=true]{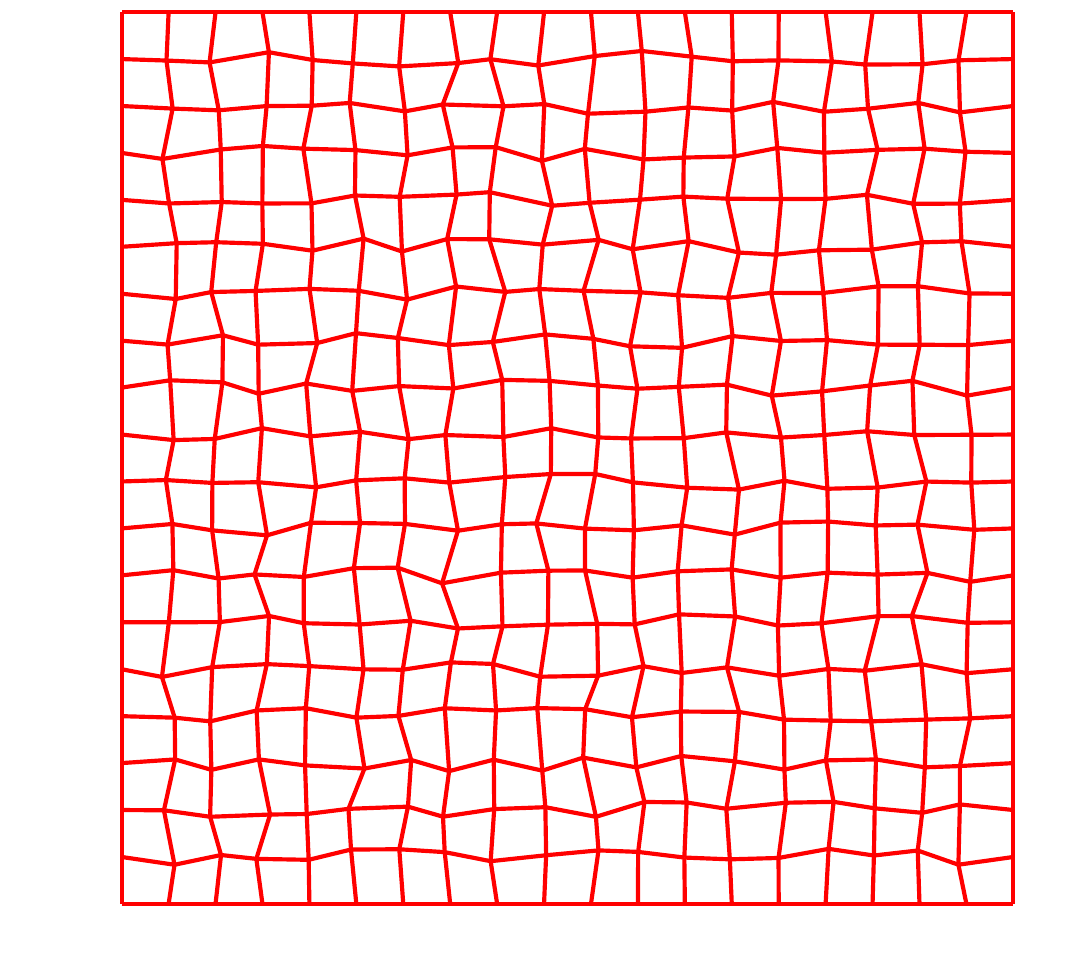}
			\caption{Average skew angle: $\alpha=\:6.1^{\circ}$}
			\label{fig:Jitter2}
		\end{subfigure}
		~
		\begin{subfigure}[b]{0.2\linewidth}
			\centering
			\includegraphics[width=\linewidth,trim= 50mm 100mm 50mm 100mm,clip=true]{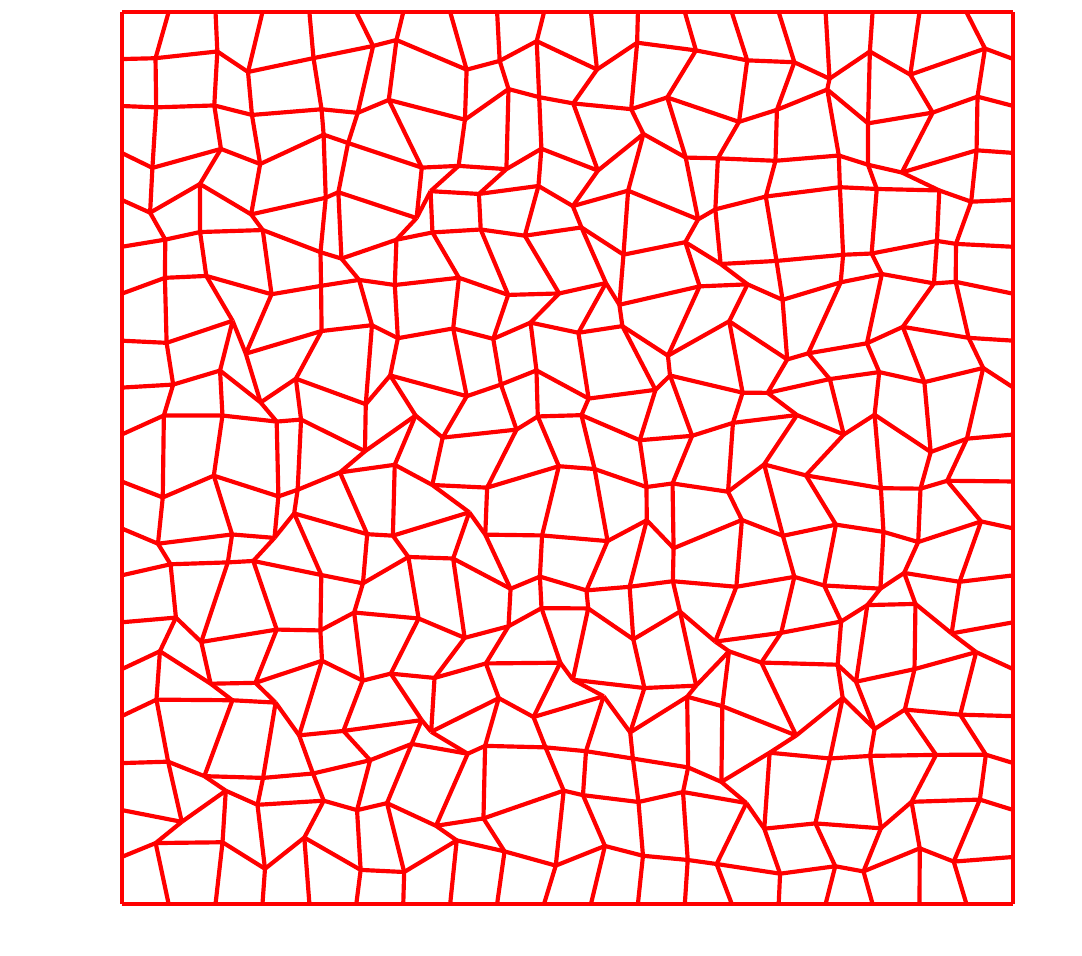}
			\caption{Average skew angle: $\alpha=\: 15.0^{\circ}$}
			\label{fig:Jitter3}
		\end{subfigure}
		\caption{$19\times19$ quadrilateral meshes showing differing degrees of node jittered mesh warp.}
		\label{fig:JitterExamples}
	\end{figure}

Few problems confronted in engineering are ever sufficiently simple that they can solved with sufficient accuracy by 1D methods, making extension to higher dimensionality crucial. Two dimensions also allows for a greater range of geometrical deformations to occur, even while maintaining a linear transformation of elements. Included within this, each element has a higher number of degrees of freedom, revealing a potential mechanism for  inarrucacies to enter the solutions. To evaluate the effect of higher dimensionality the isentropic convecting vortex (ICV) test case is used, as it has a known analytical solution, so numerical error can be straightforwardly calculated for the Euler equations. A mixture of mesh qualities are to be tested, so mesh quality was artificially reduced by  stochastically jittering corner nodes of a uniform grid via time seeded random numbers. The degree of jitter is controlled by a multiplying factor and the mesh quality and warp is then characterised by a skew angle. This is  defined as the mesh average absolute angle by which the element cross diagonals deviate from square (Fig.~\ref{fig:SkewAngle}) and  encompasses both the skewness and aspect ratio of a mesh. Some sample meshes are shown in Fig.~\ref{fig:Jitter1}~-~\ref{fig:Jitter3}.   
	
	A finite volume (FV) scheme (a simplified version of T-block), with the same explicit time integration as used in the FR calculation was used for performance comparison, as it is representative of a family of schemes widely used in industry. To evaluate the spatial error, the temporal error has to be minimised by use of an appropriately small time step, which in this case corresponds to $\mathrm{CFL} = 0.01$. ($\mathrm{CFL} = 0.05\: \& \: 0.005$ were also tested and the error was found to be independent of the temporal scheme at this level). By comparing the exact solution, $u$, and computed solution, $u^{\prime}$, the error, $\theta$, can be calculated and the spatial order of accuracy (OOA) can be obtained: 
	\begin{equation}
		\theta = u - u^{\prime} = \mathcal{O}(\delta^{p+1}) \label{eq:error}
	\end{equation}
	
\begin{figure}[h]
		\centering
			\includegraphics[width=0.5\linewidth]{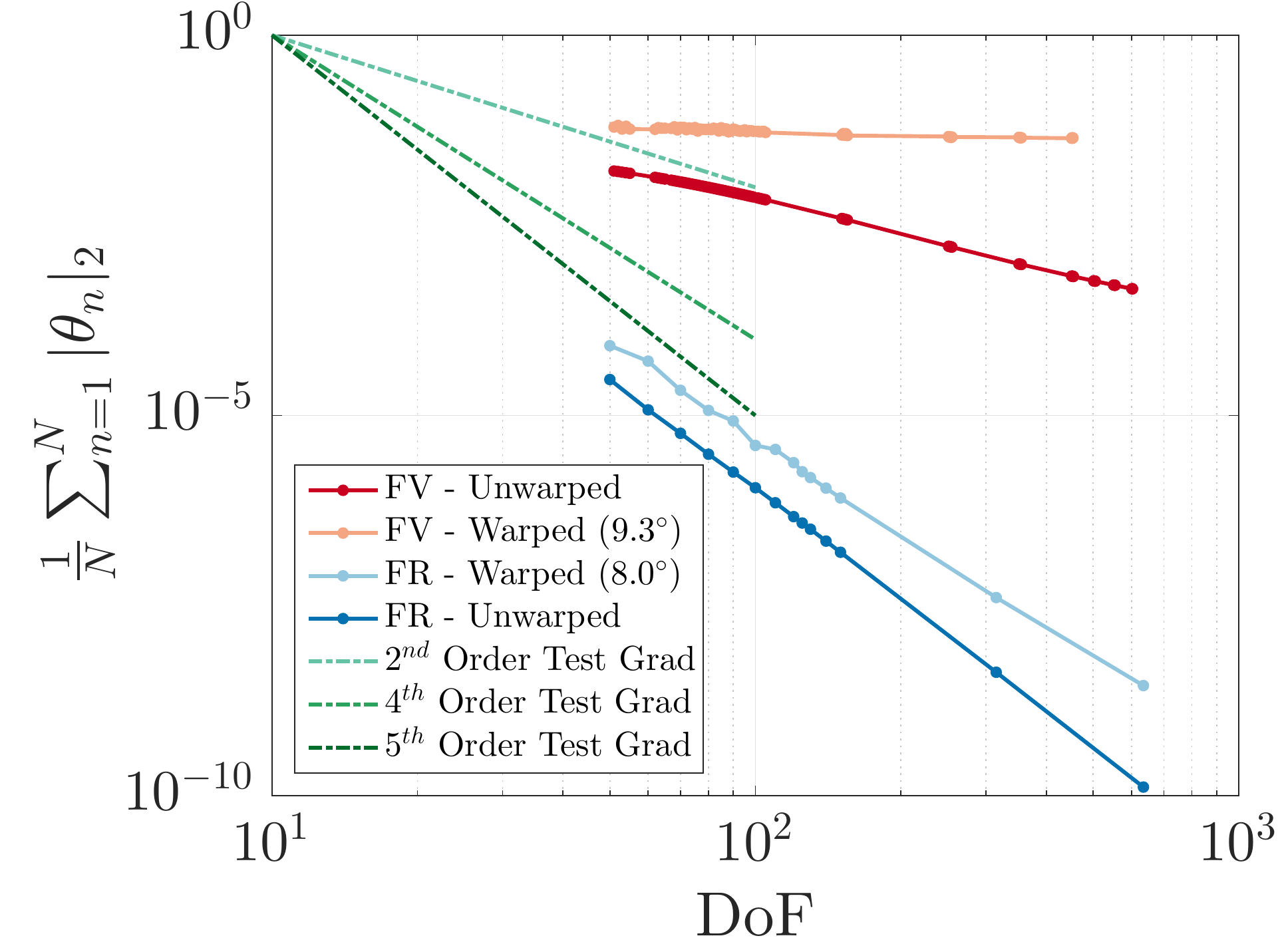}
		\captionof{figure}{Point averaged $l_2$ norm of error,$\frac{1}{N}\sum_{n=1}^N{|\theta_n|_2}$, against degrees of freedom (DoF) for $p=4$ FR and nominally second order FV scheme. $CFL = 0.01$ for $3000$ time steps.}
		\label{fig:OOA5}
	\end{figure}

	 Recovery of the spatial OOA is shown in Fig.~\ref{fig:OOA5}, \emph{i.e}, $\mathrm{OOA} = p+1$, and a comparison can be made with the $5^{th}$ order test gradient to see this. The plotted results for the FV scheme display an $\mathrm{OOA} \approx 2$, as well as a large increase in the cell averaged  $l_2$ error compared to FR for the same number of Degrees of Freedom (DoF). Also shown in Fig.~\ref{fig:OOA5} are the results of moderate mesh warping. In FR's case, a move towards $\mathrm{OOA} \approx 4$ occurs, and the FV scheme result becomes aphysical, \emph{i.e} $\mathrm{OOA} = 0$.  
	
\begin{figure}
		\centering
		\begin{subfigure}[b]{0.425\linewidth}
			\centering
			\includegraphics[width=\linewidth,trim= 32mm 92mm 41mm 98mm,clip=true]{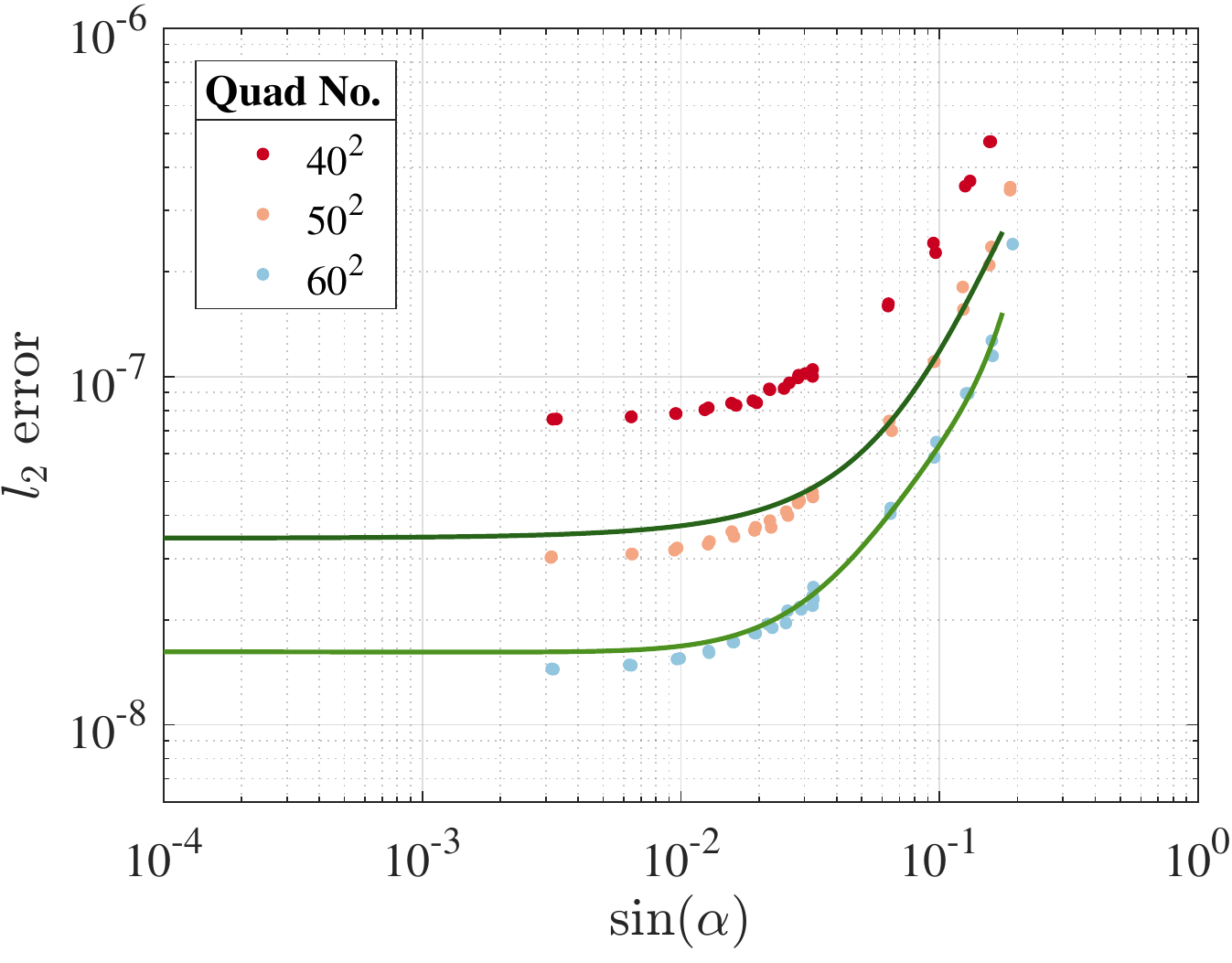}
			\caption{Magnitude of mesh warp against error for various mesh densities using an FR scheme $p=4$ ($\bullet$). Also plotted are the predicted error values for $4^{th}$ order recovery ($-$).\\ \quad}
			\label{fig:JitteredOrderFR}
		\end{subfigure}
		~
		\begin{subfigure}[b]{0.425\linewidth}
			\centering
			\includegraphics[width=\linewidth,trim= 32mm 92mm 41mm 98mm,clip=true]{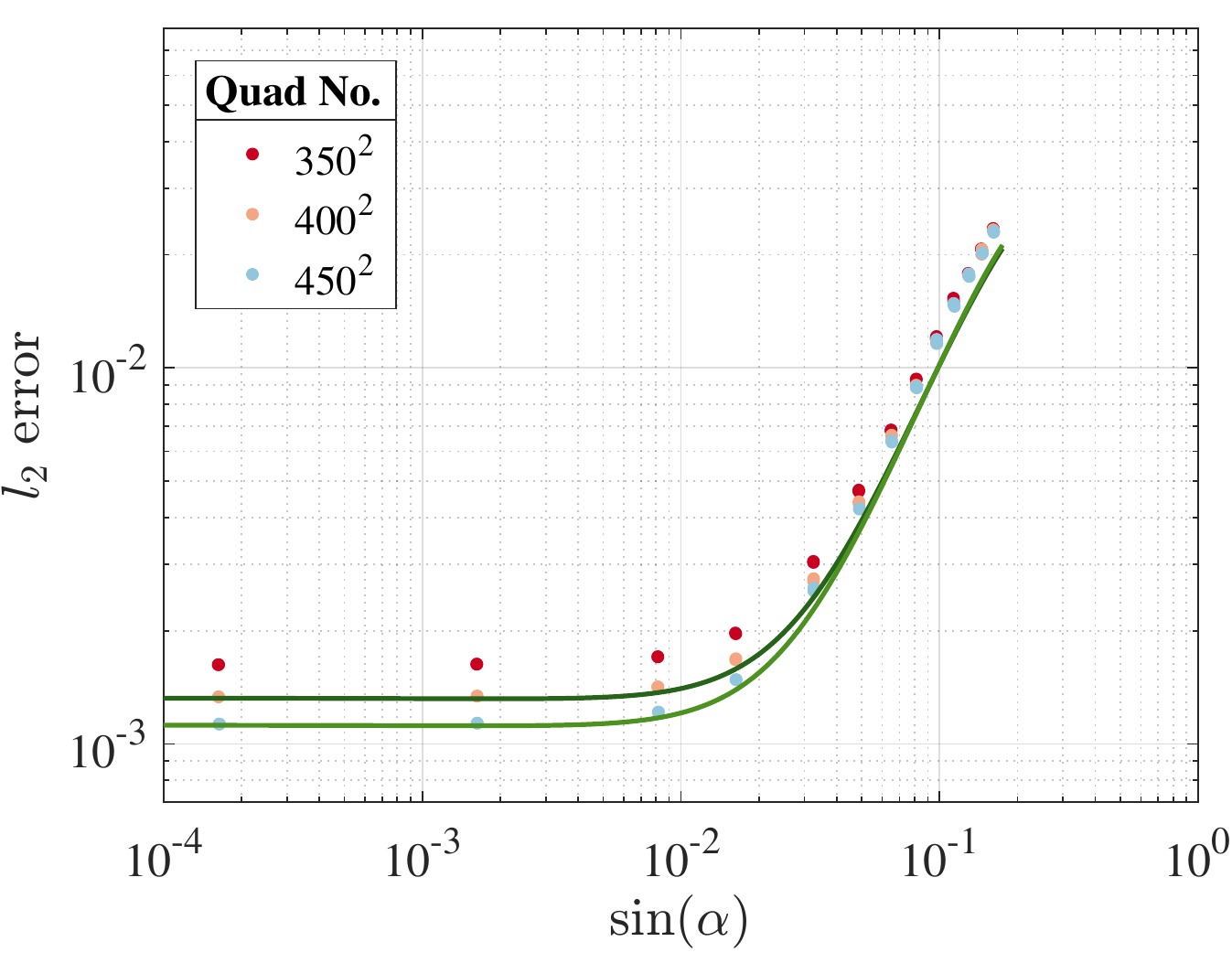}
			\caption{Magnitude of warp against error for various mesh densities using an FV scheme  ($\bullet$). Also plotted is the predicted error based on the order recovered from Fig.~\ref{fig:OOA5}.(-). The DoF are chosen to be comparable to the DoF of the FR tests.}
			\label{fig:JitteredOrderFV}
		\end{subfigure}
		\caption{}
		\label{fig:JitterOrder}
	\end{figure}

	A more detailed investigation into the deterioration of the spatial order is performed via variation of the degree of node jittering. Tests were run with $\mathrm{CFL} = 0.01$ for 500 time steps (again, results were found to be independent of CFL number at this level). The results of numerical tests and a predictive procedure are shown in Fig.~\ref{fig:JitterOrder}. The predictive procedure uses the error data from a mesh of fewer degrees of freedom together with the desired OOA to make a prediction of the error at a higher quantity of degrees of freedom.
	
	Figure~\ref{fig:JitterOrder} shows that, consistent with Fig.~\ref{fig:OOA5}, at low skew angles, FR's $\mathrm{OOA}$ is unaffected by mesh quality: $\mathrm{OOA} \rightarrow p+1$ as $\alpha \rightarrow 0^{\circ}$. Figure~\ref{fig:JitteredOrderFR} then confirms that on poorer quality meshes the $\mathrm{OOA}$ drops by approximately one, here shown as a drop from $\mathrm{OOA} =5$ to  $\mathrm{OOA} \approx 4$, this is shown by comparison to the $4^{\mathrm{th}}$ order predicted error line. This warp induced change occurs at $\alpha \approx 1.5^{\circ}$, and as mesh skewness increases to its maximum, FR is still able to give an accurate solution. By comparison the FV scheme (Fig.~\ref{fig:JitteredOrderFV}) undergoes warp induced error change at $\alpha \approx 0.55^{\circ}$, transitioning from $\mathrm{OOA} \approx 2$ to $\mathrm{OOA} = 0$. Although the FV scheme is much simpler, this shows that the performance of such schemes widely used in industry  can be rapidly eroded by a lack of mesh quality.
	
	The root cause of the loss of order accuracy of FR on warped grids is not fully revealed by these tests. The Jacobian that maps between the physical and computational domain in this case provides an exact mapping due to the node jittering providing a linear transformation. Hence, error is not introduced to the convective velocity through the Jacobian. The additional error has two potential sources. Firstly the linear component of FR could introduce error via ill-conditioning of the projection to the functional space when waves are advected at an angle. Secondly, Jameson~\etal~\cite{Jameson2012} derived the aliasing error for 1D non-linear problems in the FR framework, yet in higher dimensions, on arbitrary grids, cross multiplication of projection terms will be present. Therefore, as the mesh become skewed aliasing has the scope to introduce larger quantities of error. Both of these topics are quite expansive and are left for further investigation.
		
	As a final aside to illustrate the increase in performance that FR offers compared to a typical FV method, the ICV test case on a uniform mesh was used. For this test the number of degrees of freedom was varied such that the grid averaged $l_2$ error was comparable - this was done with a $\mathrm{CFL} = 0.01$ for $100$ time steps. (A small variation of the CFL number was made and the results were found to be invariant with CFL number). The test was carried out on a single core of an Intel\textsuperscript{\textregistered} Xeon\textsuperscript{\textregistered} L5630 which was otherwise idle. The results of testing are shown in Table~\ref{tab:FR_FV_Comp1}, and show that FR requires $\sim 2.6$ orders of magnitude less wall time for the same error in 2D as the FV method used here.   
		
	\begin{figure}
		\centering
		\captionof{table}{Comparable errors in a 2D ICV test for FR and FV schemes.}
		\label{tab:FR_FV_Comp1}
		\begin{tabular}{|l|c |c|}
			\hline
			& Flux Reconstruction & Finite Volume \\ \hline
			Point Averaged Error & $4.7\times10^{-5}$ & $4.01\times10^{-5}$ \\ \hline
			Wall Time            & $1.6\:s$           & $639\:s$            \\ \hline
			Cells                & $64,\:p=4$         & $4,000,000$         \\ \hline 
		\end{tabular}
	\end{figure}

\section{Conclusions}
	The use of FR  on warped meshes is important for the likely future applications of the scheme. It has been shown that FR is more resilient to distorted meshes than some FD and FV families of schemes which are currently widely used in industry. A detailed look at the PPW of FR on stretched grids with varying order, as well as study of the FR stability criterion, shows that, depending on geometry, the order of the scheme can be varied to increased performance. In particular, within a given cell the correction function order can be directionally varied to increase wave resolving ability. It is also shown from the linear advection equation that the CFL limit with non-regular grids is dependent on the dominant wave direction, with contracting grids providing a stabilising effect. This is  a feature of FR that will impact boundary layer meshes. It was proposed that the ill-conditioning of the functional projection of FR causes the degradation in accuracy for warped meshes. A more complete study of this is left as future work. Lastly, in some runtime comparisons, FR was found to require $\sim 2.6$ orders of magnitude less wall time for the same error in two dimensional test, compared to a widely used FV method.
\section*{Acknowledgements}
The support of the Engineering and Physical Sciences Research Council of the United Kingdom is gratefully acknowledged.

\bibliographystyle{aiaa}
\bibliography{library}

\end{document}